\documentclass[smallextended,referee,envcountsect]{svjour3}

\usepackage{float, enumerate, graphicx, amssymb, amsmath , color, longtable, multirow, lineno,url}
\usepackage{graphicx}
\usepackage[justification=centering]{caption}
\usepackage{amsmath}

\usepackage[active]{srcltx}
\usepackage{pifont}
\usepackage{dsfont}
\usepackage{paralist} 
\usepackage{amssymb}
\usepackage{authblk}
\usepackage{bbding}

\usepackage{multirow,booktabs}
\usepackage{graphicx}
\graphicspath{{figures/}}
\usepackage{subfigure}
\usepackage{epsfig}
\usepackage{epstopdf}
\usepackage{longtable,booktabs}
\usepackage[a4paper,left=2cm,right=2cm,bottom=2.3cm,top=2.3cm]{geometry}
\usepackage{lineno}
\usepackage{cite}
\usepackage[colorlinks,linkcolor=blue,citecolor=blue, bookmarksnumbered=true, bookmarksopen=true, colorlinks ]{hyperref}

\smartqed
\usepackage{graphicx}
\journalname{  }

\begin{document}

\title{\textsc{An Improved Gradient Method  with  Approximately Optimal Stepsize Based on Conic model  for Unconstrained Optimization }  }
\author{Zexian Liu    \and  Hongwei Liu    }
\institute{$\text{\underline{}   }~~$ \\	
 Zexian Liu \Envelope 
\at State Key Laboratory of Scientific and Engineering Computing, Institute of Computational
Mathematics and Scientific/Engineering computing, AMSS, Chinese
Academy of Sciences,  Beijing, 1000190,   China. 		
\at e-mail: liuzx@lsec.ac.cc.cn, liuzexian2008@163.com
\and
 Hongwei Liu
\at School of Mathematics and Statistics, Xidian University, Xi'an, 710126, People's Republic of China
\at e-mail:hwliu@mail.xidian.edu.cn
}
%
%
\date{Received: date / Accepted: date}

\maketitle

\begin{abstract} A new type of stepsize,  which was recently   introduced by Liu and Liu (Optimization, 67(3), 427-440, 2018), is   called approximately optimal stepsize    and is quit efficient for gradient method.  Interestingly,  all gradient methods can be regarded as   gradient methods with approximately optimal stepsizes.    In this paper, based on the work (Numer. Algorithms 78(1), 21-39, 2018),  we   present  an improved gradient method with    approximately optimal stepsize based on conic model  for   unconstrained optimization. If the objective function $ f $ is not close to a quadratic on the line segment between the current and latest iterates, we construct a conic model  to generate   approximately optimal stepsize for gradient method if   the conic model can be used; otherwise, we construct some quadratic models  to generate   approximately optimal stepsizes for gradient method. The convergence   of the proposed method    is  analyzed   under  suitable conditions. Numerical comparisons with some well-known conjugate gradient software packages such as CG$ \_ $DESCENT (SIAM J. Optim.   16(1), 170-192, 2005)   and CGOPT (SIAM J. Optim.  23(1), 296-320, 2013) indicate the proposed method is  very promising.
	\end{abstract}

\section{Introduction}
\vspace {-0.2cm}
\indent{  Consider the following unconstrained optimization problem
	$$\mathop {\min }\limits_{x \in {R^n}} f(x), \eqno{(1.1)} $$
	where $f:{R^n} \to R$  is  continuously differentiable and its gradient is denoted by $ g $.

The gradient method   takes  the following form
	\[{x_{k + 1}} = {x_k} - { \alpha _k}{g_k},\eqno{(1.2)} \]
where ${\alpha _k}$  is the  stepsize and $ g_k $ is the gradient of $ f $ at $ x_k $.

 Throughout this paper, $g_k = g(x_k),\; f_k = f(x_k),\; s_{k-1} = x_k-x_{k-1},\; y_{k-1} = g_k-g_{k-1}$    and $ \left\|  \cdot  \right\| $ stands for the Euclidean norm.

It is widely accepted  that the stepsize is of great importance to the numerical performance of gradient method \cite{Asmundis2013On}. In 1847, Cauchy \cite{Cauchy1847M} presented the  steepest descent   method, where the stepsize is determined by   	$$ {\alpha _k^{SD}} = \arg \;\mathop {\min }\limits_{\alpha>0}  f({x_k} - \alpha {g_k}).$$
The  steepest descent   method    usually converges   slowly. In 1988, Barzilai and Borwein \cite{Barzilai1988Two} presented   a new   gradient method (BB method), where the stepsize is given by
$$\alpha_k^{B{B_1}} = \frac{\|s_{k-1}\|^2}{s_{k-1}^Ty_{k-1}}\;\; \;\;\;\text{or} \;\;\; \;\;\; \alpha _k^{B{B_2}} = \frac{{s_{k - 1}^T{y_{k - 1}}}}{{\parallel {y_{k - 1}}{\parallel ^2}}}.  $$
Clearly, the BB method is in essence a gradient method, but  the choice of the stepsize is different from $ \alpha_k^{SD} $.

Due to the simplicity and numerical efficiency,  the BB method has enjoyed great developments during these years. The BB method has been proved to be  globally\cite{Raydan1993On} and linearly  convergent \cite{Dai2002Rlinear} for any dimensional strictly convex quadratic functions.  In 1997,  Raydan \cite{Raydan1997The} presented a global BB method for general nonlinear unconstrained optimization by incorporating the nonmonotone line search (GLL line search) \cite{Grippo1986A}, and the numerical results in \cite{Raydan1997The} suggested that the BB method is superior to some classical conjugate gradient methods. From then on,  a number of modified BB stepsizes have been exploited  for gradient methods. Dai et al. \cite{Dai2006CBB} presented the cyclic BB method for unconstrained optimization. Using the interpolation scheme, Dai et al. \cite{Dai2002Modified} presented two modified BB stepsizes for  gradient methods.   Based on some modified secant equations, Xiao et al. \cite{Xiao2010Notes} designed  four modified BB stepsizes for gradient methods.    According to a fourth order   model and some modified secant equations,  Biglari and Solimanpur \cite{Biglari2013Scaling} presented some modified gradient methods with  modified BB  stepsizes, and the numerical results in \cite{Biglari2013Scaling}  indicated that  these modified BB methods are efficient.   Miladinovi\'c et al. \cite{Miladinovic2011Scalar} proposed a new stepsize based on the usage of both the quasi-Newton property and the Hessian
inverse approximation by an appropriate scalar matrix for gradient method.

Different from the above modified BB methods, Liu and Liu \cite{Liu2018GMAOSquad} introduced a new type of stepsize for  gradient method in 2018,  which is called approximately optimal  stepsize    and  is   quite efficient for gradient method.


 \noindent  {\textbf{Definition 1.1 }  Suppose that $ f $ is continuously differentiable, and let $ {\phi_k }(\alpha ) $ be an approximation model of $ f(x_k -\alpha g_k ) $, where $  g_k  = \nabla f\left(x_k \right)  $.  A  positive number  ${\alpha_k ^ {AOS} }$ is called the  \textbf{approximately optimal  stepsize} associated to  ${\phi_k }(\alpha )$ for gradient method,  if  ${\alpha_k ^ {AOS} }$   satisfies
 	\[{\alpha_k ^  {AOS} } = \arg \;\mathop {\min }\limits_{\alpha >0}  {\phi_k }(\alpha ).\]
 	
  The approximately optimal  stepsize is generally calculated easily  and can be applied to unconstrained optimization.  In any gradient method  for strictly convex quadratic minimization problems, the   stepsize $ \alpha_k $  can also be generated by minimizing the following quadratic approximation  model:
  \[{\phi _k}\left( \alpha  \right) = {f_k} - \alpha {\left\| {{g_k}} \right\|^2} + \frac{1}{2}{\alpha ^2}g_k^T\left( {\frac{1}{{{\alpha _k}}}I} \right){g_k},\]
  where ${\frac{1}{{{\alpha _k}}}I}  $ is an approximation to the Hessian matrix of $ f $ at $ x_k $. Then, the stepsize  $ \alpha_k $ is the approximately optimal stepsizes associated to the above-mentioned  $ {\phi _k}\left( \alpha  \right) $ for gradient method. As a result,     all gradient methods can be regarded as   gradient methods with approximately optimal stepsizes in this sense.
  
We see from the definition of approximately optimal stepsize that the numerical performance of gradient method with approximately optimal stepsize depends heavily on approximation model   $ {\phi_k }(\alpha ) $. Some gradient  methods  with approximately optimal stepsizes \cite{Liu2018GMAOScone,Liu2018Several,Liu2018GMAOSNABB,Liu2019TensorBB}  
 were later proposed for unconstrained   optimization, and the numerical results in  \cite{Liu2018GMAOScone,Liu2018Several,Liu2018GMAOSNABB,Liu2019TensorBB} suggested that these gradient methods  with approximately optimal stepsizes  are  surprisingly  efficient.

In those gradient methods with approximately optimal stepsizes \cite{Liu2018GMAOScone,Liu2018Several,Liu2018GMAOSNABB,Liu2019TensorBB}, the gradient method with approximately optimal stepsizes based on conic model \cite{Liu2018GMAOScone} has enjoyed   some attentions \cite{Snezana} due to its good nice numerical performance. In this paper,     we present an improved  gradient method with   approximately optimal optimal optimal stepsize  based on conic model for  unconstrained optimization. 
In the proposed method,  when the objective function $ f $ is not close to a quadratic function on the line segment between      $ x_{k-1} $ and $ x_k $,  a conic model is exploited  to generate  approximately optimal stepsize    if the conic model   can be used. Otherwise,  some quadratic models are constructed to derive  approximately optimal stepsizes.
We   analyze the    convergence   of the proposed method  under mild conditions. Two collect sets denoted by 80pAdr and 144pCUTEr, which are from  \cite{Andrei2008An} and \cite{Gould2001CUTEr}, respectively, are used to examine the effectiveness of the test methods. Some numerical experiments  indicate that   the proposed method   is superior to the limited memory conjugate gradient software package CG$ \_ $DESCENT (6.0) \cite{HagerZhang2013The} for   80pAdr and is comparable to CG$ \_ $DESCENT (5.0) \cite{Hager2005A} for  144pCUTEr, and performs better than CGOPT \cite{Dai2013CGOPT} for   80pAdr and is comparable to CGOPT for   144pCUTEr.

The remainder of this  paper is organized as follows. In Section 2,   we exploit some approximation models including a conic model and quadratic models  to derive   efficient   approximately optimal  stepsizes for gradient method.  In Section 3, we present an improved    gradient method with     approximately optimal  stepsize based on conic model, and analyze the global convergence    of the proposed method under some suitable conditions. In Section 4, some numerical experiments are done to examine the effectiveness of the proposed method.    Conclusions  and discussions are given in the last section.
}

\section{Derivation of  the Approximately Optimal  Stepsize  }
In the section, based on the properties of the objective function $ f $,  some approximation models including a conic model and quadratic models are exploited to generate   approximately optimal  stepsize    for    gradient method.

According to the definition of approximately optimal stepsize in Section 1, we know   that the effectiveness of   approximately optimal stepsize    will rely on the  approximation model.     We determine the approximation models based on the following observations.


Define
\[{\mu _k} = \left| {\frac{{2\left( {{f_{k - 1}} - {f_k} + g_k^T{s_{k - 1}}} \right)}}{{s_{k - 1}^T{y_{k - 1}}}} - 1} \right|.\]
According to \cite{Dai2002Modified,Yuan1995A},
$ \mu_k $ is a   quantity showing
how $ f(x)  $ is close to a quadratic on the line segment between $ x_{k-1}  $ and $ x_k $.
If the following condition \cite{Liu2018GMAOScone,Liu2019TensorBB
	} holds, namely,  $$ {\mu _k} \le {c_1} \;\;\;\;  \text{or}  \;\;\;\;  \max \left\{ {{\mu _k},{\mu _{k - 1}}} \right\} \le {c_2}, \eqno{(2.1)}  $$
where  $ c_1 $ and $ c_2   $ are   small positives  and $ c_1  < c_2 $, then $ f $ might be close to a quadratic on the line
segment between $ x_{k-1}  $ and $ x_k $.   General  iterative methods, which are often based on      quadratic model,
have been quite successful in solving practical optimization problems \cite{Han2005An}, since   quadratic model   can approximates     the objective function $ f $ well at a small neighbourhood of $ x_k $ in many cases. Consequently, if $ f $ is close to a quadratic on the line segment between $ x_{k-1}  $ and $ x_k $,    quadratic   model is preferable.   However, when $x_k  $ is far   from the minimizer,   quadratic model   might   not work very well
if the objective function  $ f $ possesses high non-linearity  \cite{Sun1996Optimization,Sun2012A}. To  address the drawback,   some conic models \cite{Sun2012A,Davidon1980Conic,Sorensen1980The} have been exploited to approximate    the objective function.    The conic functions, which    interpolate both function values and gradients at the latest two iterates,   can    fit exponential functions, penalty functions or other functions  which share with conics the property of increasing rapidly near some $ n-1 $ dimensional hyperplane in $ R^n $ \cite{Davidon1980Conic}.
All of these indicate  that,  when $ f $ is not close to a quadratic function on the line segment between $ x_{k-1} $ and $ x_k $,  conic models may serve better than    quadratic model \cite{Sorensen1980The}.

Based on the above observations, we determine   approximately optimal stepsize for gradient method in the following cases.

\textbf{Case I: Conic Model  } 



 When $ f $ is not close to a quadratic   on the line segment between $ x_{k-1} $ and $ x_k $, we consider the following conic model :
 \[{\phi _k}\left( d \right) = {f_k} + \frac{{g_k^Td}}{{1 + b_k^Td}} + \frac{{{\alpha ^2}}}{2}\frac{{d_{}^T{B_k}d}}{{{{\left( {1 + b_k^Td} \right)}^2}}},\]
 where  \[{b_k} =   \frac{{1 - {\gamma _k}}}{{{\gamma _k}g_{k - 1}^T{s_{k - 1}}}}{g_{k - 1}} ,\; {\gamma _k} = \frac{{ - g_{k - 1}^T{s_{k - 1}}}}{{{\rho _k} + {f_{k - 1}} - {f_k}}},\;\;{\rho _k} = \sqrt \Delta_k  ,\;\;\Delta_k {\rm{ = }}{\left( {{f_{k - 1}} - {f_k}} \right)^2} - \left( {g_k^T{s_{k - 1}}} \right)\left( {g_{k - 1}^T{s_{k - 1}}} \right) \]
and
$ B_k $ is generated by imposing   generalized  BFGS update formula \cite{Sun2012A} on a positive scalar matrix $ D_k  $:
\[{B_k} = {D_k} - \frac{{{D_k}{v_{k - 1}}v_{k - 1}^T{D_k}}}{{v_{k - 1}^T{D_k}{v_{k - 1}}}} + \frac{{{r_{k - 1}}r_{k - 1}^T}}{{v_{k - 1}^T{r_{k - 1}}}},\]
where  $ {r_{k-1}} = {\bar y _{k-1}}/{\gamma _k} $,  $  {v_{k-1}} = {\gamma _k}{s_{k-1}} $ and $ {\bar y _{k-1}} = {\gamma _k}{g_{k  }} - \dfrac{1}{{{\gamma _{k }}}}{g_{k-1}}$.
Here we take the scalar matrix $ D_k$ as  $ D_k = \xi_1 \dfrac{{{v_{k-1}^T}v_{k-1}}}{{v_{k-1}^T{r_{k-1}}}} I $ , where $ \xi_1 \ge 1 $. It is easy to verify that, if $ {v_{k-1}^T{r_{k-1}}} > 0 $, then $D_k  $ is symmetric positive definite and thus $B_k  $ is symmetric positive definite.  In order to improve the   numerical performance,   we restrict  $ {\gamma _k} = \max \left\{ {\min \left\{ {{\gamma _k},2} \right\},0.01} \right\} $ and the coefficient $ \frac{{1 - {\gamma _k}}}{{{\gamma _k}g_{k - 1}^T{s_{k - 1}}}}  $ of $ b_k $   as $   \max \left\{ {\min \left\{ {{ \frac{{1 - {\gamma _k}}}{{{\gamma _k}g_{k - 1}^T{s_{k - 1}}}}  },5000} \right\},-5000} \right\}             $.

By substituting $ d=-\alpha g_k $ into the above conic model $ \phi_k\left(d \right)  $, we obtain that
\[ \phi_k^1 \left( {  \alpha  } \right) = f\left( {{x_k}} \right) - \frac{{\alpha g_k^T{g_k}}}{{1 - \alpha b_k^T{g_k}}} + \frac{{{\alpha ^2}}}{2}\frac{{g_k^T{B_k}{g_k}}}{{{{\left( {1 - \alpha b_k^T{g_k}} \right)}^2}}}.\]

It is clear that  $ \alpha_k = \frac{1}{b^T_kg_k} $ is the singular point of $ \phi_k^1\left( \alpha \right) $,  ${\phi _k^1}^\prime \left( \alpha  \right) = \frac{{\alpha \left( {g_k^T{B_k}{g_k} + g_k^T{g_k}b_k^T{g_k}} \right) - g_k^T{g_k}}}{{{{\left( {1 - \alpha b_k^T{g_k}} \right)}^3}}} $ and  $ \phi_k^1\left( \alpha \right) $ is continuous differentiable in $ R\backslash \left\{ {1/b_k^T{g_k}} \right\} $.

 If $  \Delta_k > 0  $,  $  {v_{k-1}^T{r_{k-1}}} > 0 $ and  $ {g_k^T{B_k}{g_k} + \left( {g_k^T{g_k}} \right)\left( {b_k^T{g_{k }}} \right)} \ne 0 $,
 by imposing $  \dfrac{{d\phi_k^1 }}{{d\alpha }}{\rm{ = }}0 $ we obtain the unique stationary point of $\phi_k^1 \left( \alpha \right)  $:
\[\alpha^S _k  =  \frac{{g_k^T{g_k}}}{{g_k^T{B_k}{g_k} + \left( {g_k^T{g_k}} \right)\left( {b_k^T{g_{k }}} \right)}}. \eqno{(2.2)}\]

 We  analyze the properties of the stationary point $ \alpha_k^S  $ in the following  two  cases.

(1) The singular point $ \alpha_k = \dfrac{1}{b_k^Tg_k} < 0 $. If $  {g_k^T{B_k}{g_k} + \left( {g_k^T{g_k}} \right)\left( {b_k^T{g_{k  }}} \right)} < 0 $, then we know $    \alpha^S_k <  \frac{1}{b_k^Tg_k} $.  By $ \phi_k^1 \left( \alpha \right) $, it is not difficult to obtain that
\[\mathop {\lim }\limits_{\alpha  \to 1/b_k^T{g_k}} {\phi _1}\left( \alpha  \right){\rm{ =  + }}\infty \;\;\;{\rm{and}}\;\;\;\;\mathop {\lim }\limits_{\alpha  \to {\rm{ + }}\infty } {\phi _1}\left( \alpha  \right){\rm{ = }}{f_k} + \frac{{g_k^T{g_k}}}{{b_k^T{g_k}}} + \frac{{g_k^T{B_k}{g_k}}}{{2{{\left( {b_k^T{g_k}} \right)}^2}}}\;\]
and $ {\phi _k^1}^\prime \left( \alpha  \right) < 0 $ for  $ \alpha>\frac{1}{b_k^Tg_k} $. Therefore, there no exists $ \alpha^* >0 $ such that $ \alpha^* = \mathop {\min }\limits_{\alpha  > 0} \;{\phi _k^1}\left( \alpha  \right) $. Consequently, if $  {g_k^T{B_k}{g_k} + \left( {g_k^T{g_k}} \right)\left( {b_k^T{g_{k  }}} \right)} \le 0 $, then  we will switch to Case II. Here we only consider the case of  $  {g_k^T{B_k}{g_k} + \left( {g_k^T{g_k}} \right)\left( {b_k^T{g_{k  }}} \right)} > 0 $.  In the case we know that $ \alpha_k^S > 0 $.  If $ \alpha > \alpha_k^S $, then  $ {\alpha \left( {g_k^T{B_k}{g_k} + g_k^T{g_k}b_k^T{g_k}} \right) - g_k^T{g_k}} > 0 $, which together with $    {1 - \alpha b_k^T{g_k}}>0  $ implies that  $${\phi _1}^\prime \left( \alpha  \right)  > 0 $$
for    $ \alpha > \alpha_k^S $.
By $ {\phi _1}^\prime \left( 0  \right) =  -\left\|g_k \right\|^2   < 0$, the continuous differentiability   of $ \phi_k^1\left(\alpha \right) $ in $ R\backslash \left\{ {1/b_k^T{g_k}} \right\} $,  the uniqueness of the stationary point  and $ {\phi _1}^\prime \left( \alpha^S_k  \right) = 0$, we know that  ${\phi _1}^\prime \left( \alpha  \right) < 0 $ holds for $ \alpha \in  \left[ 0, \alpha_k^S \right)  $. Therefore, the stationary point $ \alpha^S_k $   satisfies
\[{\alpha^S _k} = \mathop {\min }\limits_{\alpha  > 0} \;{\phi _1}\left( \alpha  \right),\]
which means that the stationary point $ \alpha^S_k $  is the approximately optimal stepsize associated to $ \phi_k^1 \left( \alpha
    \right) $.

(2)The singular point $ \alpha_k = \frac{1}{b_k^Tg_k}  > 0 $. It is obvious that the stationary point $ \alpha^S_k  $   satisfies $0 < {\alpha^S _k} < \frac{1}{{b_k^T{g_k}}}$.  If $\alpha^S_k < \alpha < \frac{1}{b^T_kg_k}  $, we obtain that $ 1-\alpha b^T_kg_k >0 $ and $ {\alpha \left( {g_k^T{B_k}{g_k} + g_k^T{g_k}b_k^T{g_k}} \right) - g_k^T{g_k}} > 0 $, which imply that
 $${\phi _1}^\prime \left( \alpha  \right)   > 0 $$
 for   $ \alpha \in \left( \alpha^S_k ,   \frac{1}{b^T_kg_k} \right)  $.
By $ {\phi _1}^\prime \left(  0  \right) = -\left\|g_k \right\| ^2< 0 $, $ {\phi _1}^\prime \left( \alpha^S_k  \right) =0 $,  the continuous differentiability   of $ \phi_k^1\left(\alpha \right) $ in $ R\backslash \left\{ {1/b_k^T{g_k}} \right\} $ and  the uniqueness of the stationary point, we know that ${\phi _1}^\prime \left( \alpha  \right)  < 0 $ holds for   $ \alpha\in \left[   0,   \alpha^S_k\right)  $. Therefore, the stationary point $ \alpha^S_k $  is a local minimizer of $ \phi_k^1 \left( \alpha \right)  $ and \[\;{\phi _1}\left( {\alpha _k^S} \right) = {f_k}  - \frac{{{{\left( {g_k^T{g_k}} \right)}^2}}}{{2g_k^T{B_k}{g_k}}}.\]
If $ \alpha > \frac{1}{b_k^Tg_k } $,  then we have $ 1 - \alpha b_k^Tg_k < 0$,
\[\mathop {\lim }\limits_{\alpha  \to {{\left( {1/b_k^T{g_k}} \right)}^ + }} {\phi _1}\left( \alpha  \right){\rm{ =  + }}\infty \;\;\;\;\;{\rm{and}}\;\;\;\;\mathop {\lim }\limits_{\alpha  \to {\rm{ + }}\infty } {\phi _1}\left( \alpha  \right){\rm{ = }}{f_k} + \frac{{g_k^T{g_k}}}{{b_k^T{g_k}}} + \frac{{g_k^T{B_k}{g_k}}}{{2{{\left( {b_k^T{g_k}} \right)}^2}}},\;\]
which together with the fact that $ {\phi _1}^\prime \left( \alpha  \right) < 0 $ holds for    $ \alpha > \frac{1}{b_k^Tg_k} $ implies that
  \[{\phi _1}\left( {\alpha _k^S} \right) - \mathop {\lim }\limits_{\alpha  \to {\rm{ + }}\infty } {\phi _1}\left( \alpha  \right) = {\phi _1}\left( {\alpha _k^S} \right) - {f_k} - \frac{{g_k^T{g_k}}}{{b_k^T{g_k}}} - \frac{{g_k^T{B_k}{g_k}}}{{2{{\left( {b_k^T{g_k}} \right)}^2}}} =  - \frac{{{{\left( {g_k^T{g_k}} \right)}^2}}}{{2g_k^T{B_k}{g_k}}} - \frac{{g_k^T{g_k}}}{{b_k^T{g_k}}} - \frac{{g_k^T{B_k}{g_k}}}{{2{{\left( {b_k^T{g_k}} \right)}^2}}} < 0  \]
holds for $ \alpha > \frac{1}{b_k^Tg_k} $.  Therefore, the stationary point $ \alpha^S_k $  satisfies
 	\[{\alpha^S_k } = \arg \;\mathop {\min }\limits_{\alpha >0}  {\phi_k^1 }(\alpha ),\]
 which implies that    the stationary point $ \alpha^S_k $  is the approximately optimal stepsize associated to $ \phi_k^1\left(
 \alpha \right) $.

 It is observed by numerical experiments that the bound $\left[ \alpha_k^{B{B_2}},\alpha_k^{B{B_1}} \right]$ for  $ \alpha^S _k    $    is very preferable for the case of  $ s_{k-1}^Ty_{k-1}>0 $.  Therefore, if the condition (2.1) does not hold and the conditions      $$  \Delta_k > 0, \;\; {v_{k-1}^T{r_{k-1}}} > 0 \;\; \text{and } \;\; {g_k^T{B_k}{g_k} + \left( {g_k^T{g_k}} \right)\left( {b_k^T{g_{k  }}} \right)} > 0  \eqno{(2.3)}  $$
hold,  the approximately optimal stepsize is     taken as  follows:
\[\alpha _k^{AOS\;(1)} = \left\{ {\begin{array}{*{20}{c}}
	{\max \left\{ {\min \left\{ {\alpha _k^S,\alpha _k^{B{B_1}}} \right\},\alpha _k^{B{B_2}}} \right\},\;\;\;{\rm{if}}\;s_{k - 1}^T{y_{k - 1}} > 0,\;\;\;\;\;}\\
	{\alpha _k^S,\;\;\;\;\;\;\;\;\;\;\;\;\;\;\;\;\;\;\;\;\;\;\;\;\;\;\;\;\;\;\;\;\;\;\;\;\;\;\;\;\;\;\;\;\;\;\;\;\;{\rm{if}}\;s_{k - 1}^T{y_{k - 1}} \le 0.\;\;\;\;}
	\end{array}} \right. \eqno{(2.4)}\]

 \textbf{Case II: Quadratic Models}

 (i)$ s_{k-1}^T y_{k-1} > 0 $

It is generally accepted that   quadratic model will serve well if $ f $ is close to  a quadratic function on the segment between  $ x_{k-1} $ and $ x_k $.  So we do not wish to abandon   quadratic model because of the large amount of practical experience and theoretical work indicating its suitability.  If  the condition (2.1)    holds  and $ s_{k-1}^Ty_{k-1} >0 $,  or    the conditions (2.2) do  not hold and $ s_{k-1}^Ty_{k-1} >0 $,
 we   consider the following quadratic approximation model:
\[\phi_k^2 (\alpha ) = f({x_k}) - \alpha \parallel {g_k}{\parallel ^2} + \frac{1}{2}{\alpha ^2}g_k^T{  B_k}{g_k},   \]
where $  B_k $ is  a  symmetric and positive definite   approximation to the Hessian matrix. 
Taking into account the storage cost and computational cost, $   B_k $  is generated by imposing the quasi-Newton update formula on a scalar matrix. Taking  the scalar matrix as $  D_k = \xi _2\frac{{y_{k - 1}^T{y_{k - 1}}}}{{s_{k - 1}^T{y_{k - 1}}}}I  $, where $ \xi _2 \ge 1 $, and imposing the modified  BFGS update formula \cite{Zhang1999NewBFGS} on the scalar matrix $   D_k  $, we obtain \[{  B_k} = {  D_k} - \frac{{{  D_k}{s_{k - 1}}s_{k - 1}^T{  D_k}}}{{s_{k - 1}^T{  D_k}{s_{k - 1}}}}{\rm{ + }} \frac{{{\bar y_{k - 1}} ^T{\bar y_{k - 1}}}}{{s_{k - 1}^T{\bar y_{k - 1}}}},   \]
where  $ {\bar y _{k - 1}} = {y_{k - 1}} +  \frac{{\bar r_k }}{\left\| {{s_{k - 1}}} \right\|^2} s_{k-1}$ and ${\bar r_k} =  {3{{\left( {{g_k} + {g_{k - 1}}} \right)}^T}{s_{k - 1}} + 6({f_{k - 1}} - {f_k})} $.

Since there exists $ u_1  \in \left[ {0,1} \right] $ such that
\[\bar r_k  = 3\left( {s_{k - 1}^T{y_{k - 1}} - s_{k - 1}^T{\nabla ^2}f({x_{k - 1}} + {u_1}{s_{k - 1}}){s_{k - 1}}} \right),\]
in order to improve the   numerical performance   we restrict $  {\bar r_k} $ as
\[\bar r_k = \min \left\{ {\max \left\{ {\bar r_k, - { \bar \eta  }s_{k - 1}^T{y_{k - 1}}} \right\},{  \bar \eta  }s_{k - 1}^T{y_{k - 1}}} \right\},  \eqno{(2.5)}\]
where $ 0 <  \bar \eta   <  0.1 $.

It follows from (2.5) that $ {s_{k - 1}^T{\bar y_{k - 1}}}  = s_{k-1}^Ty_{k-1} + \bar r_k \ge (1-  \bar \eta   ) s_{k-1}^Ty_{k-1}   $ when $ s_{k-1}^Ty_{k-1} > 0 $, which implies the following lemma.

\textbf{Lemma 2.1 } Suppose that $ s_{k-1}^Ty_{k-1} > 0 $. Then $ {s_{k - 1}^T{\bar y_{k - 1}}} > 0 $ and $B_k $  is symmetric and  positive definite.

Imposing $ \frac{{d{\phi _k^2}}}{{d\alpha }} = 0 $, we obtain
\[\bar \alpha_k^{AOS(2)}   = \frac{{g_k^T{g_k}}}{{g_k^T{B_k}{g_k}}} = \frac{{g_k^T{g_k}}}{{\frac{{{\xi _2}{{\left\| {{y_{k - 1}}} \right\|}^2}}}{{s_{k - 1}^T{y_{k - 1}}}}\left( {{{\left\| {{g_k}} \right\|}^2} - \frac{{{{\left( {g_k^T{s_{k - 1}}} \right)}^2}}}{{{{\left\| {{s_{k - 1}}} \right\|}^2}}}} \right) + {{\left( {\frac{{g_k^T{y_{k - 1}}}}{{{{\left( {s_{k - 1}^T{{\bar y}_{k - 1}}} \right)}^2}}} + \frac{{{{\bar r}_k}g_k^T{s_{k - 1}}}}{{{{\left( {s_{k - 1}^T{{\bar y}_{k - 1}}} \right)}^2}{{\left\| {{s_{k - 1}}} \right\|}^2}}}} \right)}^2}}}. \eqno{(2.6)}\]
By $ s_{k-1}^Ty_{k-1}>0 $ and Lemma 2.1, we know that $\bar \alpha_k^{AOS(2)} $     is the approximately optimal stepsize associated to $ \phi_k^2 (\alpha) $.

 It is also observed by numerical experiments that the bound $\left[ \alpha_k^{B{B_2}},\alpha_k^{B{B_1}} \right]$ for  $ \alpha _k  $ in (2.6)   is very preferable. Therefore, if the condition (2.1)  holds and $s_{k - 1}^T{y_{k - 1}} > 0$,   or the conditions (2.3) do  not hold and  $s_{k - 1}^T{y_{k - 1}} > 0$, the approximately optimal stepsize is  taken as the  truncation form of   $ \bar \alpha_k^{AOS(2)}   $ :
 \[  \alpha _k^{AOS\;(2)} = \max \left\{ {\min \left\{ {\bar \alpha_k^{AOS(2)}  ,\alpha _k^{B{B_1}}} \right\},\alpha _k^{B{B_2}}} \right\}. \eqno{(2.7)}\]

 (ii)$ s_{k-1}^T y_{k-1} \le 0 $ 

It is a   challenging task  to determine a suitable stepsize for gradient method when $ s_{k-1}^Ty_{k-1}  \le 0.$ In some modified BB  methods \cite{Dai2002Modified,Xiao2010Notes},  the  stepsize  is set simply to ${\alpha _k} =10^{30}$ for the case of   $ s_{k-1}^Ty_{k-1} \le 0 $. It is too simple to consume  expensive computational cost for searching a suitable stepsize   for gradient method.

In \cite{Liu2018GMAOSNABB}, Liu et al.  proposed a  simple and efficient strategy for choosing the stepsize for the case of  $ s_{k-1}^Ty_{k-1}  \le 0$: $ \alpha_k=\delta \alpha_{k-1} $, where $ \delta >0 $.
Liu and Liu \cite{Liu2018Several} designed an approximation model to generate  approximately optimal stepsize. Liu and Liu \cite{Liu2018GMAOScone} designed two approximation models to generate two approximately optimal stepsizes, and the numerical results in \cite{Liu2018GMAOScone}
showed that these approximately optimal stepsize are efficient. We    take the stepsize \cite{Liu2018GMAOScone} for gradient method, which is described  here for completeness.

If the condition (2.1)  holds and $ s_{k-1}^Ty_{k-1} \le 0 $, or    the conditions (2.3) do  not hold and $ s_{k-1}^Ty_{k-1} \le 0 $,
we design other approximation models to derive   approximately optimal stepsizes. Suppose for the moment that $f$ is twice continuously differentiable, the second order Taylor  expansion is
\[f({x_k} - \alpha {g_k}) = f({x_k}) - \alpha g_k^T{g_k} + \frac{1}{2}{\alpha ^2}g_k^T{\nabla ^2}f({x_k}){g_k} + o({\alpha ^2}).\]
For a very small $\tau_k  >0$, we  have that
\[{\nabla ^2}f({x_k} ){g_k} \approx - \frac{{g({x_k} - \tau_k {g_k}) - g({x_k})}}{\tau_k },\;\]
\noindent{ which gives a  new approximation model
	\[\phi_k^3  \left( \alpha\right)  = f({x_k}) - \alpha g_k^T{g_k} + \frac{1}{2}{\alpha ^2}|g_k^T(g({x_k} - \tau_k {g_k}) - g({x_k}))/{\tau_k} |.\]	
	\noindent{If ${g_k^T(g({x_k} - \tau_k {g_k}) - g({x_k}))/{\tau_k }}\ne0$, then by    imposing $\dfrac{{d\phi_k^3}}{{d\alpha }} = 0$ and  the coefficient of $ \alpha^2 $ in $ \phi_k^3 \left( \alpha \right) $, we  obtain the approximately optimal stepsize associated to $ \phi_k^3  \left( \alpha\right) $:
		\[\alpha^{AOS(3)}_k  = \frac{{  g_k^T{g_k}}}{|{g_k^T(g({x_k} - \tau_k {g_k}) - g({x_k}))/{\tau_k }|}}. \eqno{(2.8)}\]
		
\noindent{When  ${g_k^T(g({x_k} - \tau_k {g_k}) - g({x_k}))/{\tau_k }}=0$, similar to   \cite{Liu2018GMAOSNABB}, the stepsize  $ \alpha_k $    $ \alpha_k $  is computed by
	$$\alpha_k=\delta \alpha_{k-1}, \eqno{(2.9)}$$
	where $\delta >0 $.
	
\indent{To obtain  the  stepsize $ \alpha_k $ in (2.8),   it  has the cost of an extra gradient evaluation, which may result in great computational cost if the gradient evaluation is evoked frequently. To reduce the  computational cost,  we turn   to consider $g_{k-1}. $ Since
			\[s_{k - 1}^T{y_{k - 1}} =  - {\alpha _{k - 1}}g_{k - 1}^T({g_k} - {g_{k - 1}}) = {\alpha _{k - 1}}(\parallel {g_{k - 1}}{\parallel ^2} - g_{k - 1}^T{g_k}) \le 0,\]
			we have that
			\[\parallel {g_{k - 1}}{\parallel ^2} \le g_{k - 1}^T{g_k},\]
			which implies
			\[\frac{{\parallel {g_{k - 1}}\parallel }}{{\parallel {g_k}\parallel }} \le 1.\]
				
\noindent{If $\dfrac{{\parallel {g_{k-1}}{\parallel ^2}}}{{\parallel {g_{k }}{\parallel ^2}}} \ge \xi_3 $,   where $ \xi_3 > 0 $ is close to 1,   we know that $g_k$  and   $g_{k-1}$ will incline  to be collinear and  $\parallel g_k \parallel$  and $\parallel g_{k-1} \parallel$ are approximately     equal. In the case,    we use $  g_{k - 1}^T{\nabla ^2}f({x_k}){g_{k - 1}} $ to approximate $ g_k^T{\nabla ^2}f({x_k}){g_k}  $, and then use $  \frac{{\left| {{{\left( {(g({x_k} + {\alpha _{k - 1}}{g_{k - 1}}) - g({x_k})} \right)}^T}{g_{k - 1}}} \right|}}{{{\alpha _{k - 1}}}}=\frac{{\left| {s_{k - 1}^T{y_{k - 1}}} \right|}}{{\alpha _{k - 1}^2}} $ to estimate   $   g_{k - 1}^T{\nabla ^2}f({x_k}){g_{k - 1}} $, which imply   a new approximation model:
\[{\phi_k^4}\left( \alpha  \right)= f({x_k}) - \alpha \parallel {g_k}{\parallel ^2} + \frac{1}{2}{\alpha ^2}\frac{{\left| {s_{k - 1}^T{y_{k - 1}}} \right|}}{{\alpha _{k - 1}^2}}.\]
\noindent{If  $s_{k - 1}^Ty_{k - 1} \ne 0$, by   imposing $ \frac{{d{\phi _k^4}}}{{d\alpha }} = 0 $  and  the coefficient of $ \alpha^2 $ in $ \phi_k^4 \left( \alpha \right) $, we also obtain the approximately optimal stepsize associated to  $ \phi_k^4 \left( \alpha \right)    $:	
						\[\alpha^{AOS(4)}_k  = \frac{{\parallel {g_k}{\parallel ^2}}}{|{s_{k - 1}^T{y_{k - 1}}}|}{\alpha _{k - 1}^2}.\eqno{(2.10)}\]
						
\noindent{As for the case of  $s_{k - 1}^Ty_{k - 1}=0$, the  stepsize is also computed by (2.9).

Therefore, if the condition (2.1)  holds and $ s_{k-1}^Ty_{k-1} \le 0 $, or    the conditions (2.3) do  not hold and $ s_{k-1}^Ty_{k-1} \le 0 $, the  stepsize is determined by
\[{\alpha _k} = \left\{ {\begin{array}{*{20}{c}}
	{	 \dfrac{{g_k^T{g_k}}}{{|g_k^T(g({x_k} - {\tau _k}{g_k}) - g({x_k}))/{\tau _k}|}}\;,\;\;\;\;\text{if}\;\;\dfrac{{\parallel {g_k}{\parallel ^2}}}{{\parallel {g_{k - 1}}{\parallel ^2}}} < {\xi _3}\;\;\text{and}\;\;g_k^T(g({x_k} - {\tau _k}{g_k}) - g({x_k}))/{\tau _k} \ne 0,\;}\\		
		{\dfrac{{\parallel {g_k}{\parallel ^2}}}{{|s_{k - 1}^T{y_{k - 1}}|}}\alpha _{k - 1}^2,\;\;\;\;\;\;\;\;\;\;\;\;\;\;\;\;\;\;\;\;\;\;\;\;\;\;\;\;\;\text{if}\;\;\dfrac{{\parallel {g_k}{\parallel ^2}}}{{\parallel {g_{k - 1}}{\parallel ^2}}} \ge {\xi _3}\;\;\;\;\;\text{and}\;\;\;\;s_{k - 1}^T{y_{k - 1}} \ne 0,\;\;\;\;\;\;\;\;\;\;\;\;\;\;\;\;\;\;\;\;\;\;\;\;\;\;\;}\\	
	{\delta {\alpha _{k - 1}},\;\;\;\;\;\;\;\;\;\;\;\;\;\;\;\;\;\;\;\;\;\;\;\;\;\;\;\;\;\;\;\;\;\;\;\;\;\;\;\;\;\;\;\begin{array}{*{20}{c}}
		{\text{otherwise},\;\;}&{}&{}
		\end{array}\;\;\;\;\;\;\;\;\;\;\;\;\;\;\;\;\;\;
		\;\;\;\;\;\;\;\;\;\;\;\;\;\;\;\;\;\;\;\;\;\;\;\;\;\;\;\;\;\;\;\;\;\;\;\;\;\;\;\;\;\;\;\;\;\;\;\;\;\;\;}
	\end{array}} \right.\eqno{(2.11)}\]
where $\delta >0$.							
						

\section{Gradient Method with  Approximately Optimal  Stepsize Based on Conic Model }

In the section, we present an improved gradient method with approximately optimal stepsize based on conic model (we call it  GM$\_$AOS (cone) for short) for unconstrained optimization.  Though GLL line search \cite{Grippo1986A} was firstly incorporated into    the BB method  \cite{Raydan1997The},  it is observed by numerical experiments  that for modified BB  methods the nonmonotone line search (Zhang-Hager line search) proposed by Zhang and Hager  \cite{Zhang2004A} is preferable. Usually, the strategy  (3.2) for a nonmonotone line search \cite{Birgin2000Nonmonotone} is  used to accelerate   the convergence rate. Therefore, we adopt Zhang-Hager line search with the strategy  (3.2)   in GM$\_$AOS  (cone). Motivated by SMCG$ \_ $BB \cite{LiuSMCGBB}, at the  first iteration  we   choose the initial stepsize $ \alpha _0^0 $ as
	 \[\alpha _0^0 = \left\{ {\begin{array}{*{20}{c}}
	 	{\begin{array}{*{20}{c}}
	 		{1,}&{\;\;\;\;\;\;\;\;\;\;\;\;\;\;\;\;\;\;\;\;\;\;\;\;\;\;\;\;\;\;\;\;\;\;\;\;\;\;\;\;\;\;\;\;\;\;\;\;\;\;\;\;\;\;\;\;\;{\rm{if}}\left| f \right| \le {{10}^{ - 30}}\;{\rm{and}}\;\;{{\left\| {{x_0}} \right\|}_\infty } \le {{10}^{ - 30}},\;\;\;\;\;\;\;\;\;\;\;\;\;\;\;\;\;\;\;\;\;}
	 		\end{array}}\\
	 	{\begin{array}{*{20}{c}}
	 		{\begin{array}{*{20}{c}}
	 			{\begin{array}{*{20}{c}}
	 				{2\left| f \right|/\left\| {{g_0}} \right\|,}&{\;\;\;\;\;\;\;\;\;\;\;\;\;\;\;\;\;\;\;\;\;\;\;\;\;\;\;\;\;\;\;\;\;\;\;\;\;\;\;\;\;\;{\rm{if}}\left| f \right| > {{10}^{ - 30}}\;{\rm{and}}\;\;{{\left\| {{x_0}} \right\|}_\infty } \le {{10}^{ - 30}},\;\;\;\;\;\;\;\;\;\;\;\;\;\;\;\;\;\;\;\;\;}
	 				\end{array}}\\
	 			{\begin{array}{*{20}{c}}
	 				{\min \left\{ {1,{{\left\| {{x_0}} \right\|}_\infty }/{{\left\| {{g_0}} \right\|}_\infty }} \right\},\;}&{\;\;\;\;\;\;\;\;\;\;\;\;\;\;\;\;\;\;\;\;\;\;{\rm{if}}\;{{\left\| {{g_0}} \right\|}_\infty } < {{10}^7}\;{\rm{and}}\;\;{{\left\| {{x_0}} \right\|}_\infty } > {{10}^{ - 30}},\;\;\;\;\;\;\;\;\;\;\;\;\;\;\;\;\;\;\;\;}
	 				\end{array}}
	 			\end{array}}\\
	 		{\begin{array}{*{20}{c}}
	 			{\min \left\{ {1,\max \left\{ {1,{{\left\| {{x_0}} \right\|}_\infty }} \right\}/{{\left\| {{g_0}} \right\|}_\infty }} \right\},\;\;\;\;\;\;\;\;}&{{\rm{if}}\;{{\left\| {{g_0}} \right\|}_\infty } \ge {{10}^7}\;{\rm{and}}\;\;{{\left\| {{x_0}} \right\|}_\infty } > {{10}^{ - 30}}.\;\;\;\;\;\;\;\;\;\;\;\;\;\;\;\;\;\;\;}
	 			\end{array}}
	 		\end{array}}
	 	\end{array}} \right.\eqno{(3.1)}\]
	  Now we describe GM$\_$AOS (cone)    in detail.

\noindent{\textbf{Algorithm 1 GM$\_$AOS (cone)        } \\
	\textbf{Step 0} Initialization.Given a starting point  ${x_0} \in {R^n} $, constants   $\; \varepsilon > 0, \;   \lambda_{\min}, \lambda_{\max},  \; \eta_{\min}, \; \eta_{\max}, \; \sigma,    \; \delta, \;  \xi_1, \; \xi_2, \; \xi_3, \;$
	
\ \ \	$  {\bar \eta  },  \;   c_1$ and  $ c_2  $.   Set $  {Q_0}= 1, \;   {C_0} = {f_0} $ and $   k := 0. $ \\
	\textbf{Step 1} If $||{g_k}||_\infty \le \varepsilon $, then stop. \\
	\textbf{Step 2} Compute the initial stepsize for Zhang-Hager line search.
	
	\textbf{Step 2.1}  If $k = 0$ , then compute $ \alpha_0^0  $ by (3.1) and set   $ \alpha   =  \alpha_0^0 $,  go to Step 3.
	
	\textbf{Step 2.2} If   the   condition  (2.1) does not  hold and the conditions (2.3) hold, then compute  $\alpha _k $ by (2.4).  Set
	
\ \ \ \	 $ \alpha _k^0 = \max \left\{ {\min \left\{ {{\alpha _k},{\lambda _{\max }}} \right\},{\lambda _{\min }}} \right\} $ and $\alpha = \alpha_k^0,$  and  go to Step 3.
	
	\textbf{Step 2.3} If   $   s_{k-1}^T{y_{k-1}} > 0   $, then compute $\alpha _k $ by (2.7); otherwise compute $ \alpha_k $ by  (2.11).   Set
	
\ \ \ \ 	 $ \alpha _k^0 = \max \left\{ {\min \left\{ {{\alpha _k},{\lambda _{\max }}} \right\},{\lambda _{\min }}} \right\} $	and $\alpha = \alpha_k^0,$ and  go to Step 3.\\
	\textbf{Step 3} Zhang-Hager line search. If
	\[f({x_k} - \alpha {g_k}) \le {C_k} - \sigma \alpha \parallel {g_k}{\parallel ^2},\;\;\]
	\ \ \ \ 	\ \ \ \ 	then go to Step 4. Otherwise, update  $\alpha $  by \cite{Birgin2000Nonmonotone}
	$$\alpha  = \left\{ {\begin{array}{*{20}{c}}
		{\bar \alpha  ,\;\;\;\;\;\;\;\text{if} \; \alpha  > 0.1\alpha_k^0\;\text{and}\;\bar \alpha   \in [0.1\alpha_k^0,0.9\alpha ],\;}\\
		{0.5\alpha ,\;\;\;\text{otherwise},\;\;\;\;\;\;\;\;\;\;\;\;\;\;\;\;\;\;\;\;\;\;\;\;\;\;\;\;\;\;\;\;\;\;\;\;\;\;}
		\end{array}} \right. \eqno{(3.2)} $$
	\indent {\ \ \ where $\bar \alpha $  is the trial stepsize obtained by a quadratic interpolation at  ${x_k}$ and ${x_k} - \alpha {g_k}$, go to Step 3.\\
		\textbf{Step 4}  Choose $\eta_k \in [\eta_{\min},\eta_{\max}]$ and update $Q_{k+1}, C_{k+1}$ by the following ways:  $$Q_{k+1}= \eta_kQ_k+1, C_{k+1}=(\eta_kQ_kC_k+f(x_{k+1}))/Q_{k+1}.\eqno{(3.3)}$$
		
		\noindent {\textbf{ Step 5} Set ${\alpha _k} = \alpha $,  ${x_{k + 1}} = {x_k} - {\alpha _k}{g_k}$,  $k := k + 1$ and go to Step 1.\newline
			
  In what follows, we analyze  the convergence and the convergence rate of  GM$ \_ $AOS  (cone). Our convergence result utilizes the following assumptions : \\
  A1. $f$ is continuously differentiable on $R^n$. \\
  A2. $f$ is bounded below on $R^n$. \\
  A3. The gradient $g$  is  Lipschitz continuous on $R^n$, namely, there exists $L>0$ such that
  $$\parallel g(x) - g(y)\parallel  \le L\parallel x - y\parallel ,\;\;\forall x,y \in {R^n}.$$
 \indent{Since $ d_k =-g_k $, we have $\|d_k\|=\|g_k\|$ and $g_k^Td_k =-\|g_k\|^2$. Therefore, by Theorem 2.2 of \cite{Zhang2004A} we can easily obtain the following  theorem which shows  that  GM$ \_ $AOS (cone)  is globally convergent.
 	
 \indent{\textbf{Theorem 3.1} Suppose that assumption A1, A2 and A3 hold. Let $\{ {x_k}\} $ be the sequence generated by  GM$ \_ $AOS (cone). Then
  \[\liminf_{k \to \infty}    \left\| g_k \right\|   = 0.\]
  Furthermore, if  $\eta_{\max} < 1 ,$ then \[\mathop {\lim }\limits_{k \to \infty } \left\|  g_k \right\|   = 0.\]
  Hence, every convergent subsequence of the $\{ {x_k}\} $ approaches a stationary point $x^\ast$ .

Similar to the above theorem, by Theorem 3.1 of \cite{Zhang2004A},  we also obtain the following theorem which implies  the   R-linear convergence of GM$ \_ $AOS (cone).

 	\noindent  {\textbf{Theorem 3.2}  Suppose that A1  and A3 hold, $f$ is strongly convex with unique minimizer $x^\ast$ and $\eta_{\max} < 1$. Then there exists $ \zeta \in (0, 1)$ such that  \[f({x_k}) - f({x^*}) \le {\varsigma ^k}\left( {f({x_0}) - f({x^*})} \right),\] for each $k \ge 0$.

\section{Numerical Experiments}

In the section,    some numerical experiments are conducted to check the numerical performance of GM$ \_ $AOS (cone).  Two groups of collect  sets are used, and their  names    are described in Table 1 and Table 2, respectively. The first group of collect  sets   denoted  by 80pAndr   includes 80 test functions  mainly  from \cite{Andrei2008An}, and their     expressions    and   Fortran  codes   can be found in   Andrei's website:  \url{http://camo.ici.ro/neculai/AHYBRIDM }. The dimension of each test function in  80pro$ \_ $Andrei  is set to 10000 and the initial points are default. The second group of collect  sets   denoted by  144pCUTEr   includes  145 test functions from CUTEr  library \cite{Gould2001CUTEr},   which can be found in \url{http://users.clas.ufl.edu/hager/papers/CG/results6.0.txt}. It is noted that the 144 test functions from CUTEr  library \cite{Gould2001CUTEr} are indeed used  to test, as the default initial point is the optimal point in the test function ``FLETCBV2'',  so the second group of collect sets is denoted by 144pCUTEr. And the initial points and dimensions of the test functions from 144pCUTEr are default.

 The BB method, the SBB4 method \cite{Biglari2013Scaling},  CGOPT \cite{Dai2013CGOPT} and CG$ \_ $DESCENT\cite{Hager2005A}  are chosen to be compared with   GM$ \_ $AOS (cone). All test methods are implemented by C language.  The C code of GM$ \_ $AOS (cone) and some    numerical results  can be downloaded from the website: \url{http://web.xidian.edu.cn/xdliuhongwei/en/paper.html}. The   codes of CGOPT and  CG$ \_ $DESCENT  can be downloaded  from   \url{http://coa.amss.ac.cn/wordpress/?page_id=21} and \url{http://users.clas.ufl.edu/hager/papers/Software}, respectively.

 In the numerical experiments,   GM$ \_ $AOS (cone) uses the following parameters:   $ \varepsilon  = {10^{ - 6}},\;  \delta = 10, \; \sigma ={10^{ - 4}},\; {\lambda _{\min }} = {10^{ - 30}}, \; {\lambda _{\max }} = {10^{30}},  \; \eta_k = 1 $, $ \; \bar \eta  = 5.0/3  \times 10^{-5} $, $\; \tau_k=\min\{0.1\alpha_{k-1},0.01\}  $, $\xi_1 =2.15, \; \xi_2 =1.07, \;  \xi_3 =0.9, \; c_1=10^{-8} \;$ and $  c_2=0.07  $.   The BB method and the SBB4 method adopt the same line search as  GM$ \_ $AOS (cone). All  the gradient   methods are  stopped  if   $\parallel {g_{k}}{\parallel _\infty } \le {10^{ - 6}}$ is satisfied, the number of iterations  exceeds $140000$, or the number of function  evaluations exceeds $50000$.  CG$\_$DESCENT  and CGOPT  are terminated if   $\parallel {g_{k}}{\parallel _\infty } \le {10^{ - 6}}$ is satisfied or the number of iterations  exceeds $140000$, and use all default   parameter values in their codes but the above stopping conditions.

	The numerical experiments with 80pro$ \_ $Andrei     are running on Microsoft Visual Studio 2012, which is installed  in  Windows 7 in  a PC with 3.20 GHz CPU processor, 4 GB RAM memory, while  the numerical experiments  with 144pCUTEr  are running on Ubuntu 10.04 LTS  fixed  in  a VMware Workstation 10.0, which is installed  in  Windows 7 in  the same PC. The performance profiles introduced by Dolan and Mor\'e \cite{Dolan2002Benchmarking} are used to display the performance of these  methods, respectively.     In the following figures, $ N_{iter} $,  $ N_{f} $, $ N_{g} $ and $ T_{CPU} $ represent the performance profiles in term of the number of iterations, the number of function evaluations, the number of gradient evaluations and CPU time (s), respectively.
	
	The numerical experiments  are divided into four groups.
	
 	\begin{figure*}[htbp]
 		\begin{minipage}[t]{0.5\linewidth}
 			\includegraphics[width=8cm,height=5.5cm]{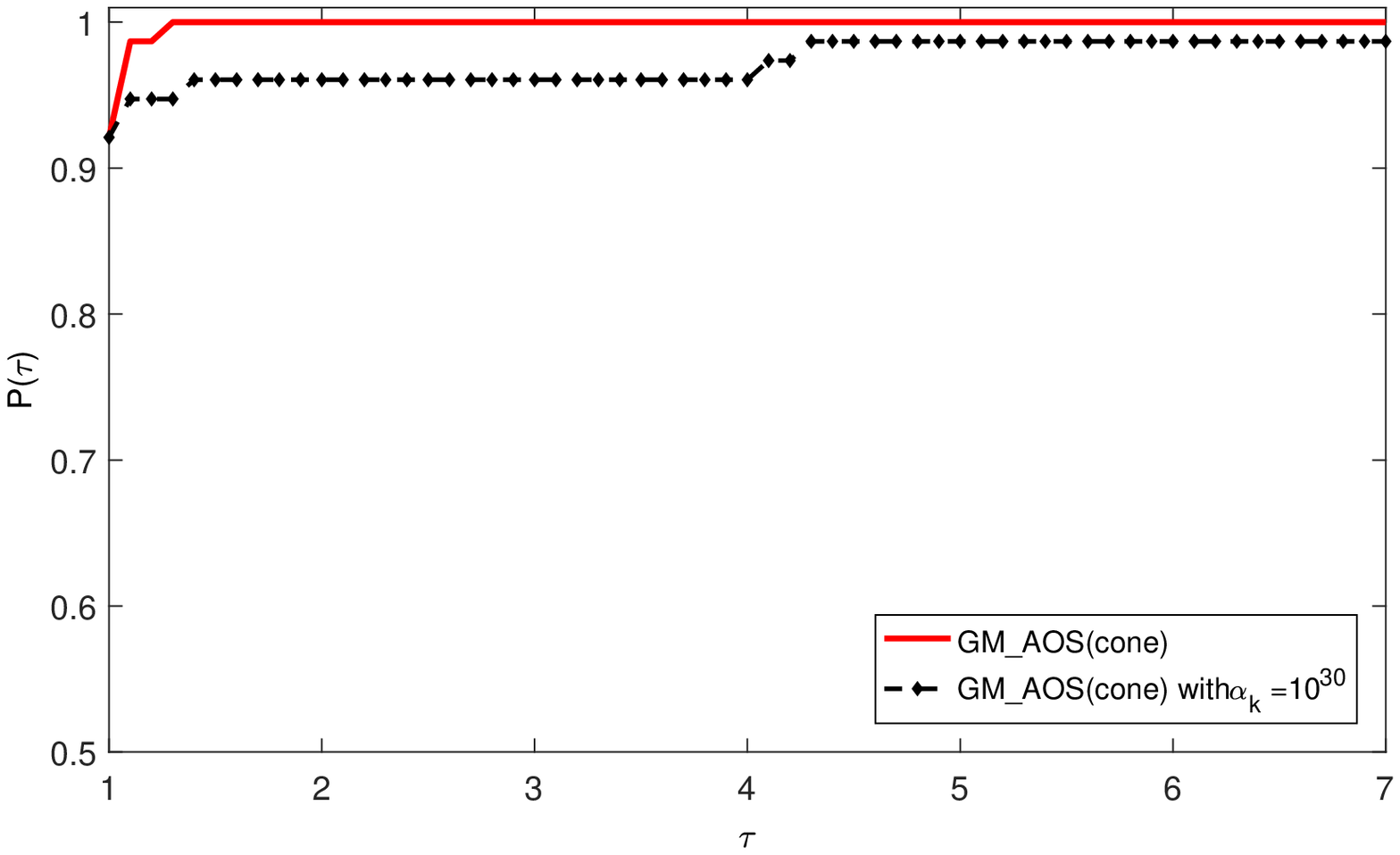}
 			\caption{    $ N_{iter} $  (80pAndr)}
 			\label{fig:5-1}
 		\end{minipage}%
 		\begin{minipage}[t]{0.5\linewidth}
 			\includegraphics[width=8cm,height=5.5cm]{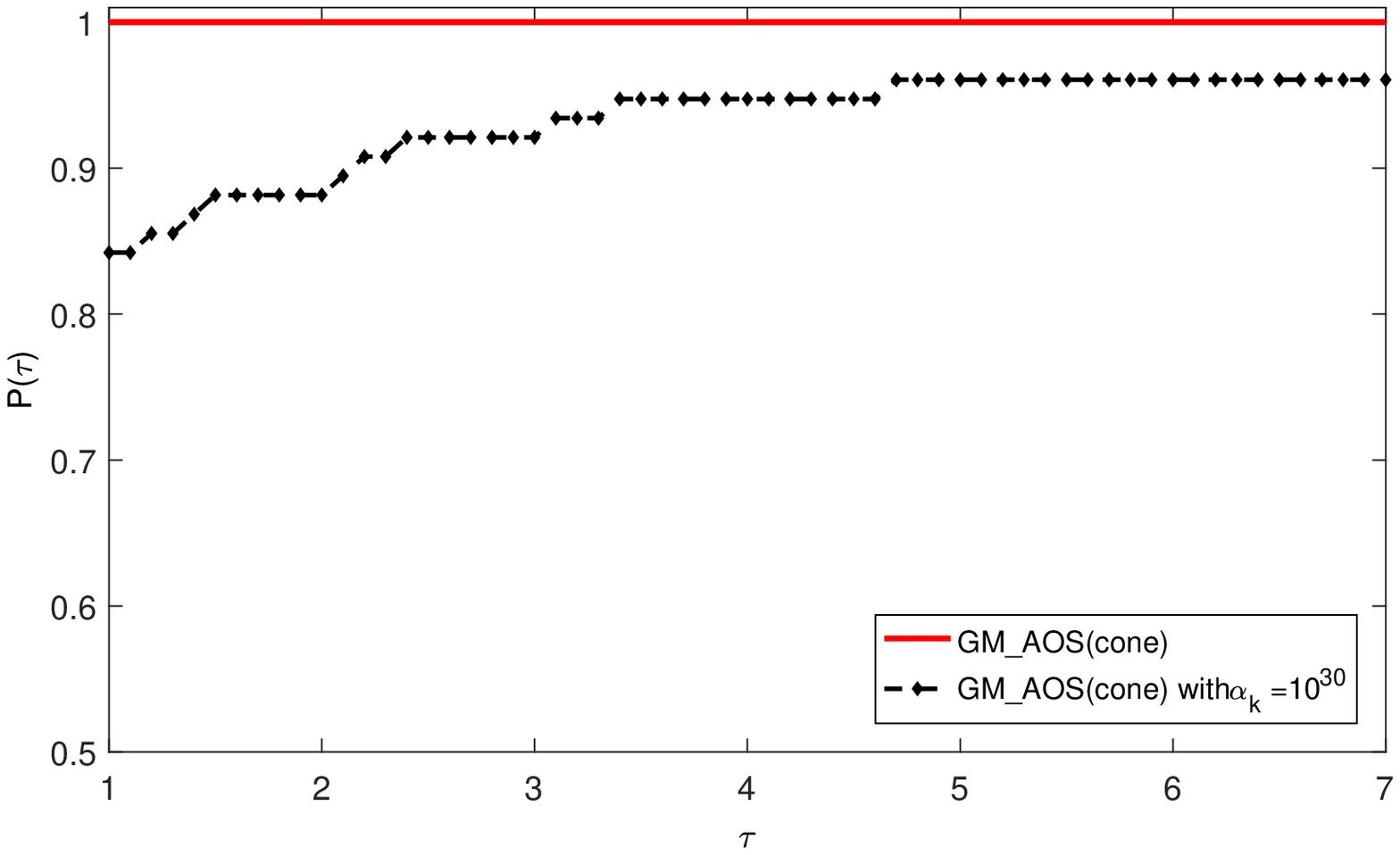}
 			\caption{   $ N_f $  (80pAndr)}	\label{fig:5-2}
 		\end{minipage}%
 	\end{figure*}
 	\begin{figure*}[htbp]
 		\begin{minipage}[t]{0.5\linewidth}
 			\includegraphics[width=8cm,height=5.5cm]{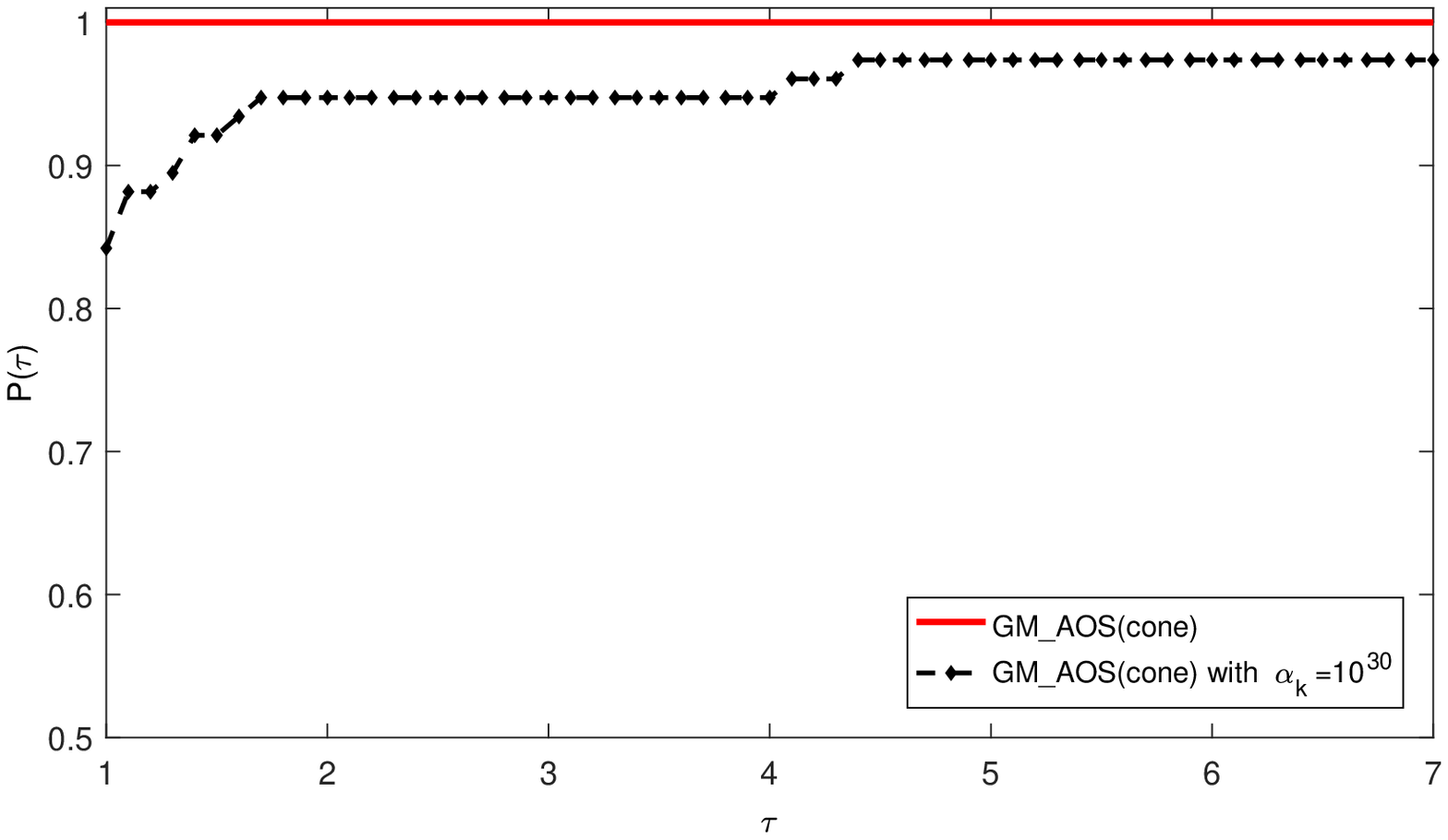}
 			\caption{   $ N_g $  (80pAndr)}
 			\label{fig:5-3}
 		\end{minipage}%
 		\begin{minipage}[t]{0.5\linewidth}
 			\includegraphics[width=8cm,height=5.5cm]{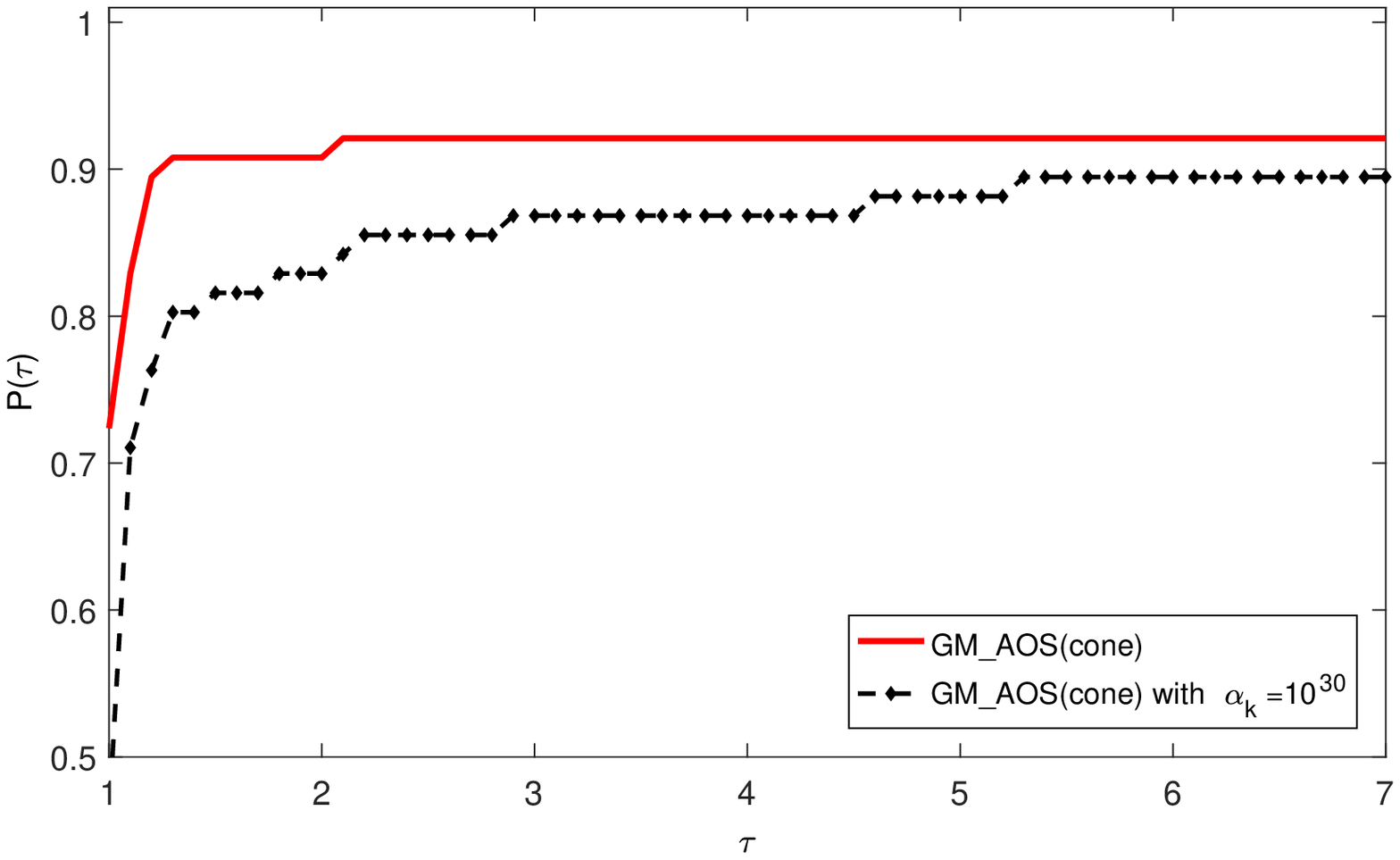}
 			\caption{   $ T_{CPU} $  (80pAndr)}	\label{fig:5-4}
 		\end{minipage}%
 	\end{figure*}

 \indent{ In the first group of numerical experiments, we  use the collect set 80pAndr to   examine the effectiveness of the stepsize (2.11).   In Figs. \ref{fig:5-1}-\ref{fig:5-4}, ``GM$ \_ $AOS (cone) with $ \alpha_k=10^{30} $'' stands for the variant of GM$ \_ $AOS (cone), which is different from  GM$ \_ $AOS (cone) only in   that (2.11) is replaced by $ \alpha_k=10^{30} $  in the Step 2.3 of GM$ \_ $AOS (cone).  In numerical experiments, GM$ \_ $AOS (cone) successfully solves all 80 problems, while its variant successfully solves 76 problems.   As shown in Fig. \ref{fig:5-1},  GM$ \_ $AOS (cone)   performs slightly better than its variant  in term of the    number  of iterations.  We can observe from Figs. \ref{fig:5-2}-\ref{fig:5-3} that   GM$ \_ $AOS (cone) requires  much  less    function   evaluations and less gradient evaluations than  its variant since the stepsize (2.11) is used. In Fig. \ref{fig:5-4}, we      see    that      GM$ \_ $AOS (cone)  is much  faster than its variant. It indicates that the stepsize (2.11) is very efficient.

		\begin{figure*}[htbp]
			\begin{minipage}[t]{0.5\linewidth}
				\includegraphics[width=8cm,height=5.5cm]{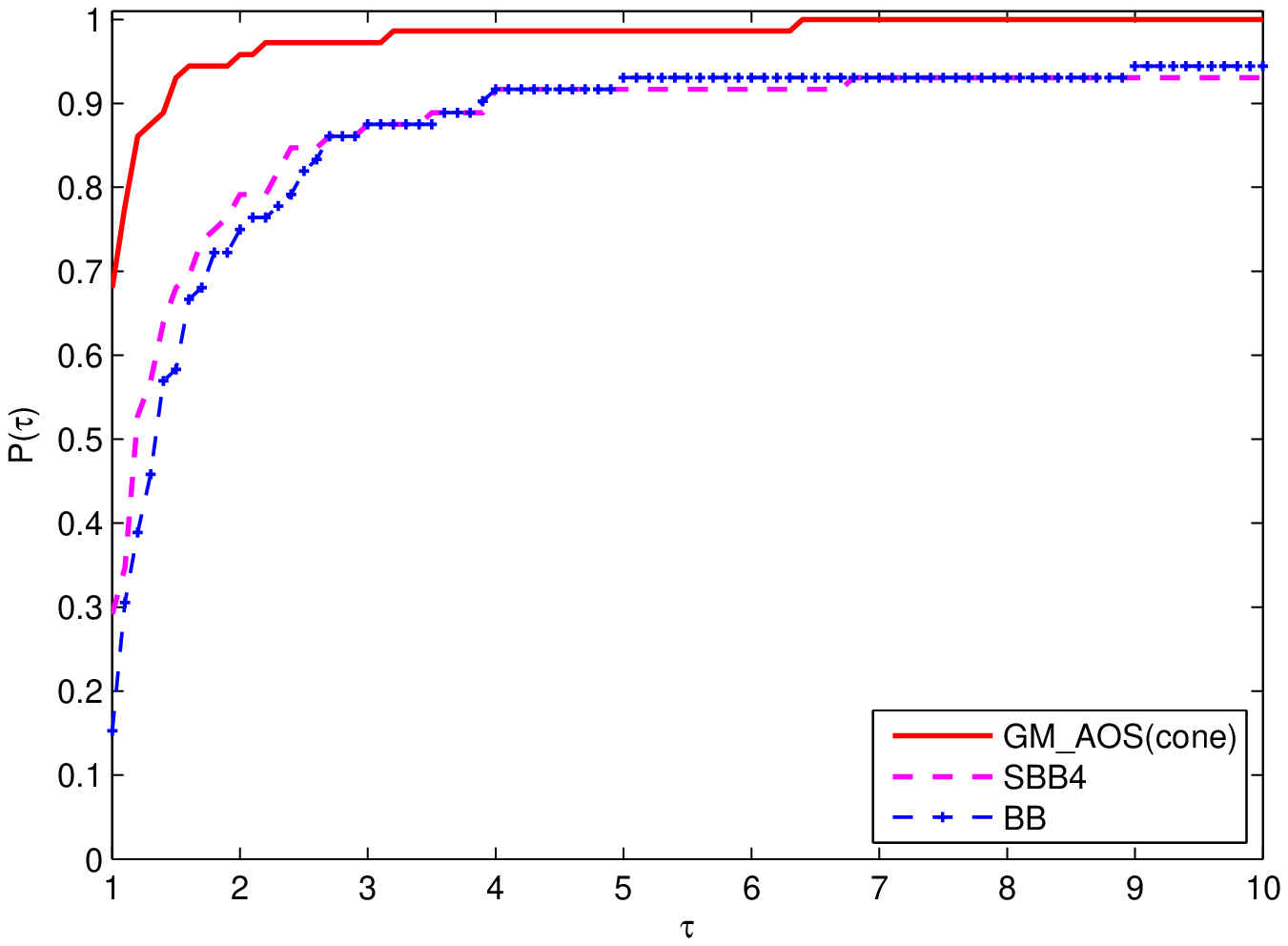}
				\caption{   $ N_{iter} $  (80pAndr)} 	\label{fig:5-5}
			\end{minipage}%
			\begin{minipage}[t]{0.5\linewidth}
				\includegraphics[width=8cm,height=5.5cm]{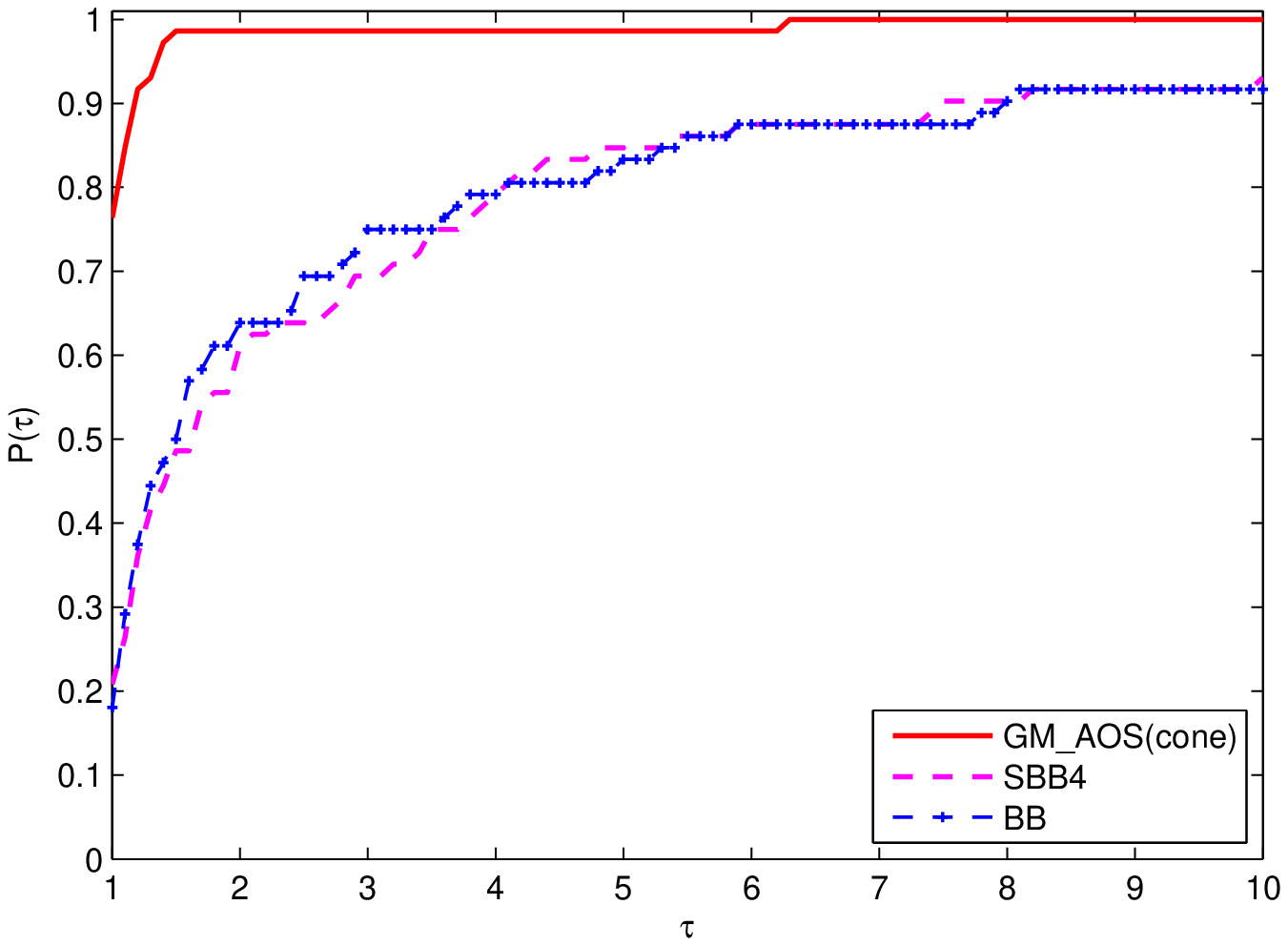}
				\caption{   $ N_{f} $  (80pAndr)}\label{fig:5-6}
			\end{minipage}%
		\end{figure*}
		\begin{figure*}[htbp]
			\begin{minipage}[t]{0.5\linewidth}
				\includegraphics[width=8cm,height=5.5cm]{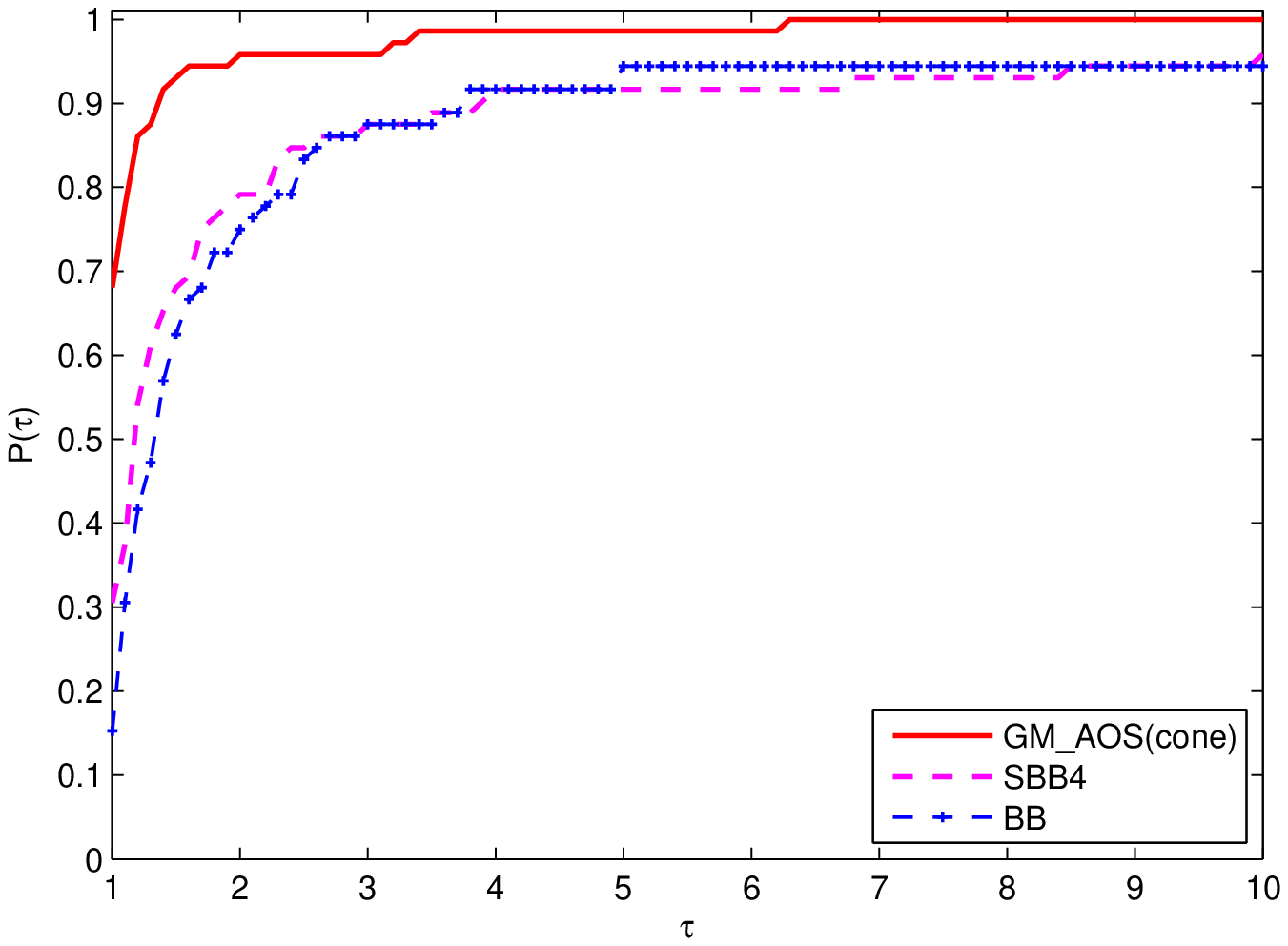}
				\caption{   $ N_{g} $  (80pAndr)}
				\label{fig:5-7}
			\end{minipage}%
			\begin{minipage}[t]{0.5\linewidth}
				\includegraphics[width=8cm,height=5.5cm]{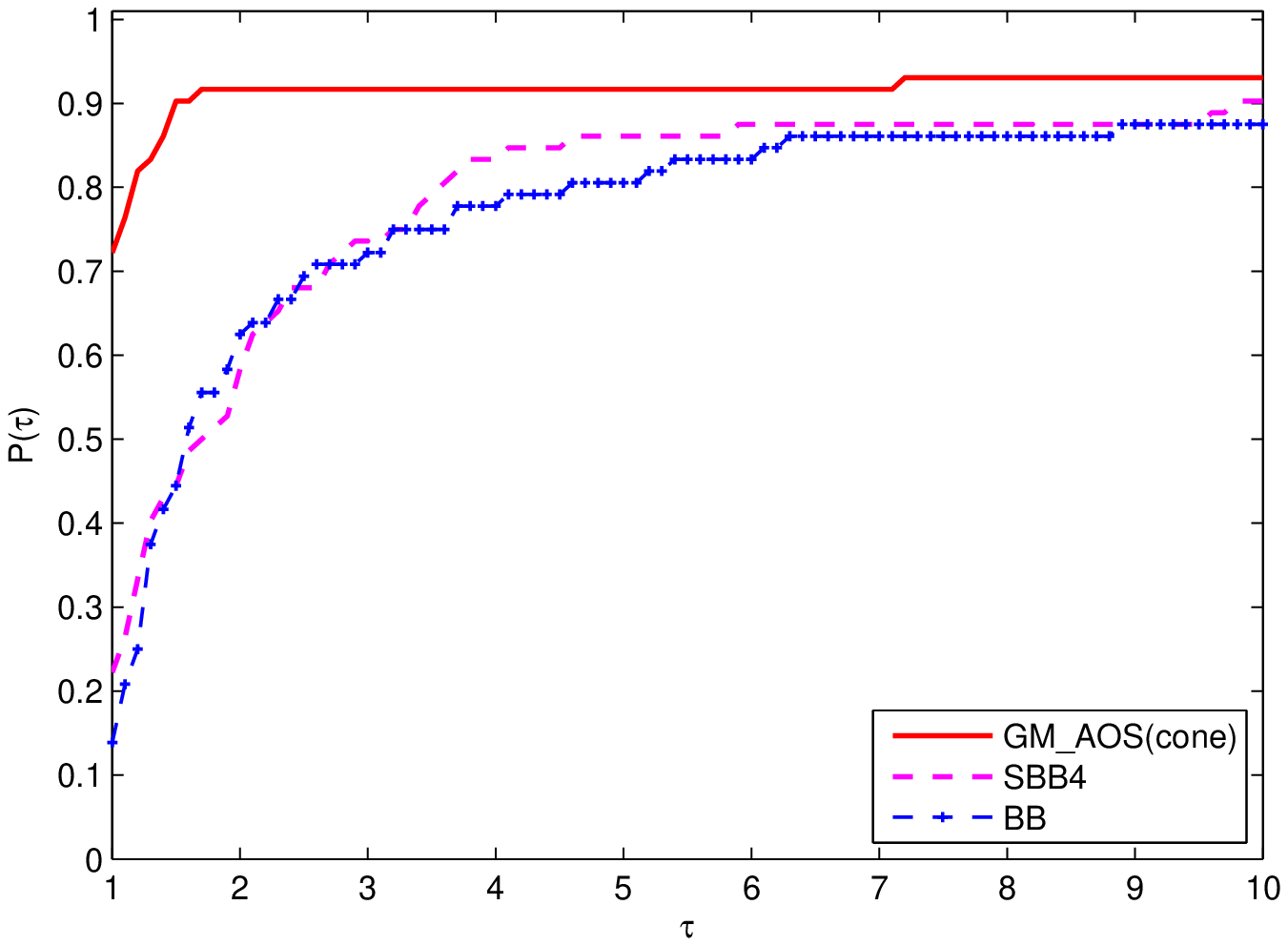}
				\caption{   $ T_{CPU} $  (80pAndr)}\label{fig:5-8}
			\end{minipage}%
		\end{figure*} 	

  In the second group of numerical experiments,  we also use the collect set 80pAndr to  compare the performances of GM$ \_ $AOS (cone)  with that of the SBB4 method and the BB method.  In numerical experiments, GM$ \_ $AOS (cone) successfully solves all 80 problems, while the SBB4 method and the BB method  successfully   solve  75 and 76 problems, respectively.
   As shown in Fig. \ref{fig:5-5}, GM$ \_ $AOS (cone) outperforms  the SBB4 method and the BB method, since  GM$ \_ $AOS (cone)  successfully solves about 68$ \% $ problems with the least iterations, while the  percentages of the SBB4 method and the BB method are about 28$ \% $ and   15$ \% $, respectively.  Similar observation   can be    made in Fig. \ref{fig:5-7} for the number of gradient evaluations.  We observe from Fig. \ref{fig:5-6} that GM$ \_ $AOS (cone)
  has a very great advantage over the SBB4  method and the BB method in term of   the number of function evaluations,    since    GM$ \_ $AOS (cone)   successfully solves about  $ 77\% $ problems  with the least function   evaluations, while the percentage of the SBB4 method and the BB
  method are $21\% $ and  $18\% $, respectively.  Fig. \ref{fig:5-8} shows that GM$ \_ $AOS (cone)  is much faster than the SBB4 method and the BB method. It indicates that GM$ \_ $AOS (cone)   is superior to the SBB4 method and the BB method.

\begin{figure*}[htbp]
	\begin{minipage}[t]{0.5\linewidth}
		\includegraphics[width=8cm,height=5.5cm]{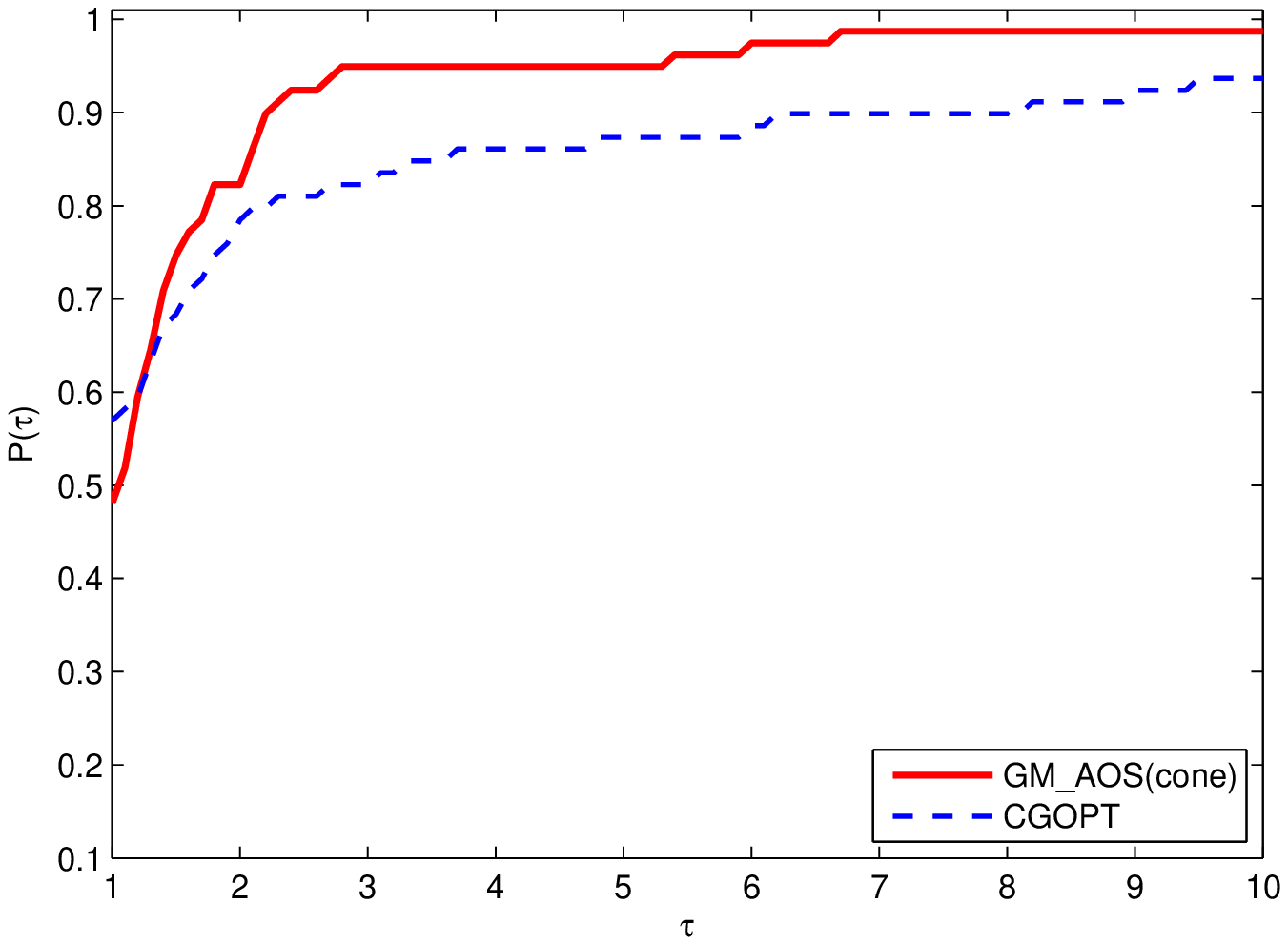}
		\caption{   $ N_{iter} $  (80pAndr)}\label{fig:5-9}
	\end{minipage}%
	\begin{minipage}[t]{0.5\linewidth}
		\includegraphics[width=8cm,height=5.5cm]{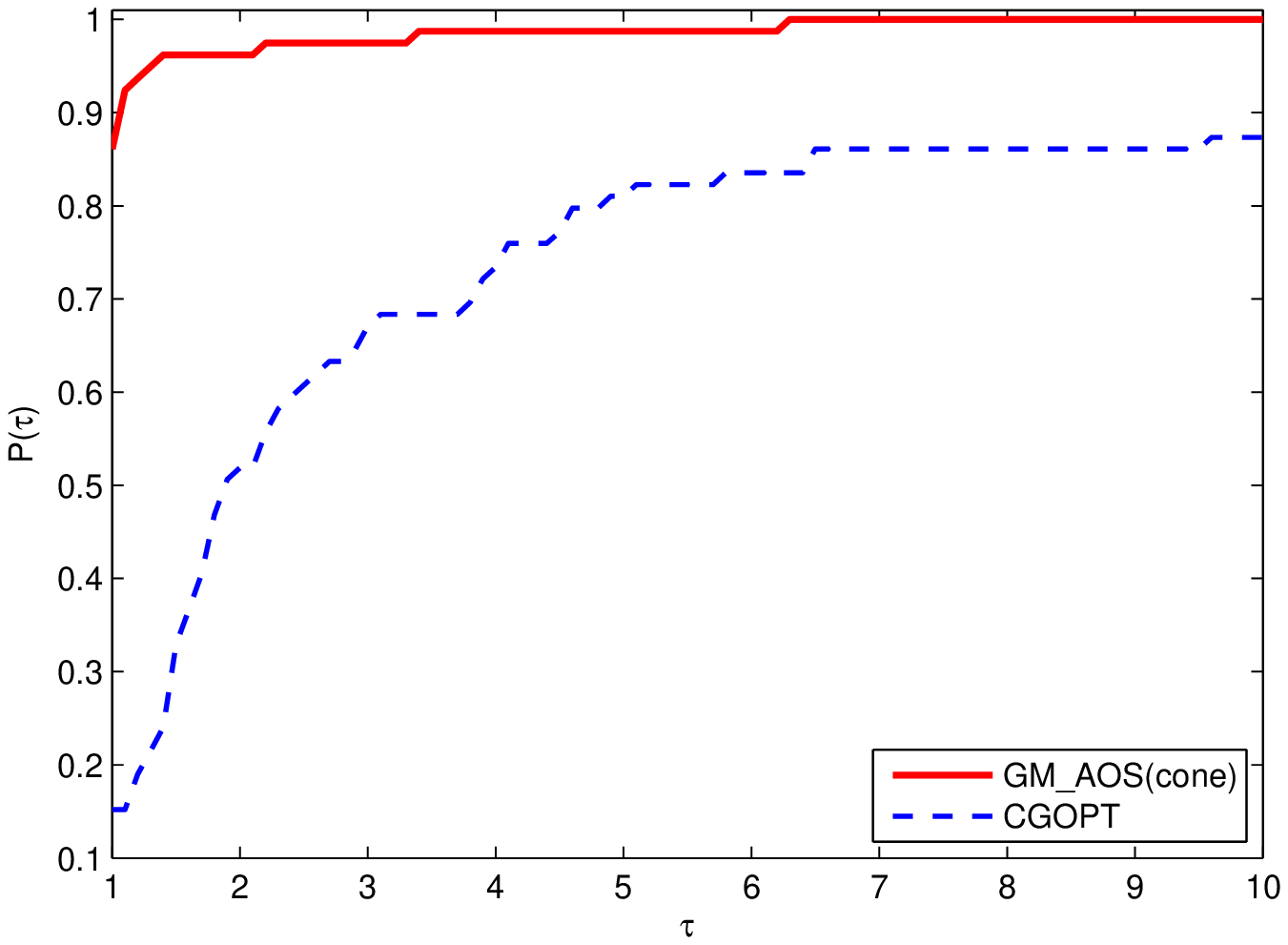}
		\caption{   $ N_{f} $  (80pAndr)}
		\label{fig:5-10}
	\end{minipage}%
\end{figure*}  	

\begin{figure*}[htbp]
	\begin{minipage}[t]{0.5\linewidth}
		\includegraphics[width=8cm,height=5.5cm]{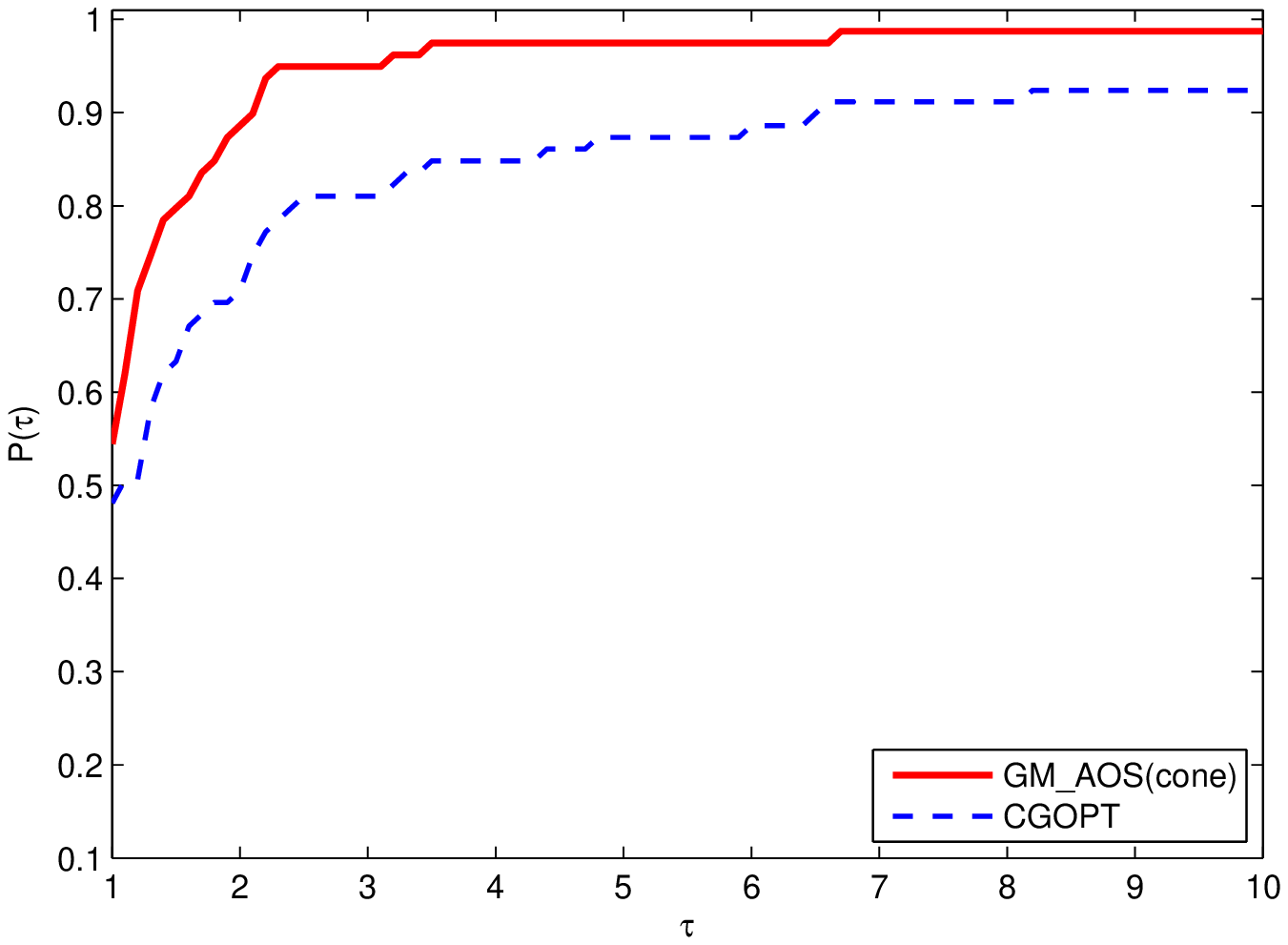}
		\caption{   $ N_{g} $  (80pAndr)}\label{fig:5-11}
	\end{minipage}%
	\begin{minipage}[t]{0.5\linewidth}
		\includegraphics[width=8cm,height=5.5cm]{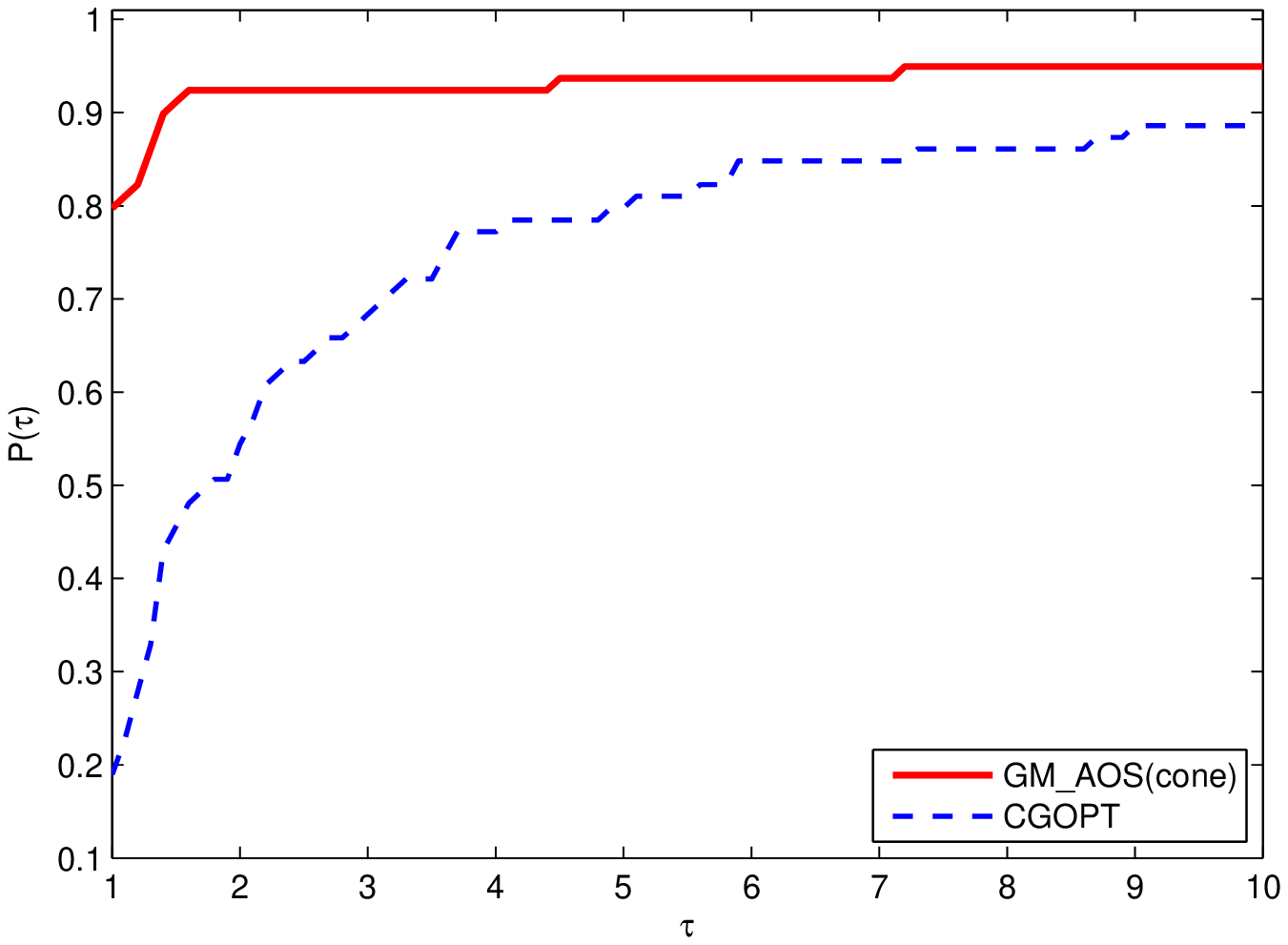}
		\caption{   $ T_{CPU} $  (80pAndr)}
		\label{fig:5-12}
	\end{minipage}%
\end{figure*}
                 \begin{figure*}[htbp]
                 	\begin{minipage}[t]{0.5\linewidth}
                 		\includegraphics[width=8cm,height=5.5cm]{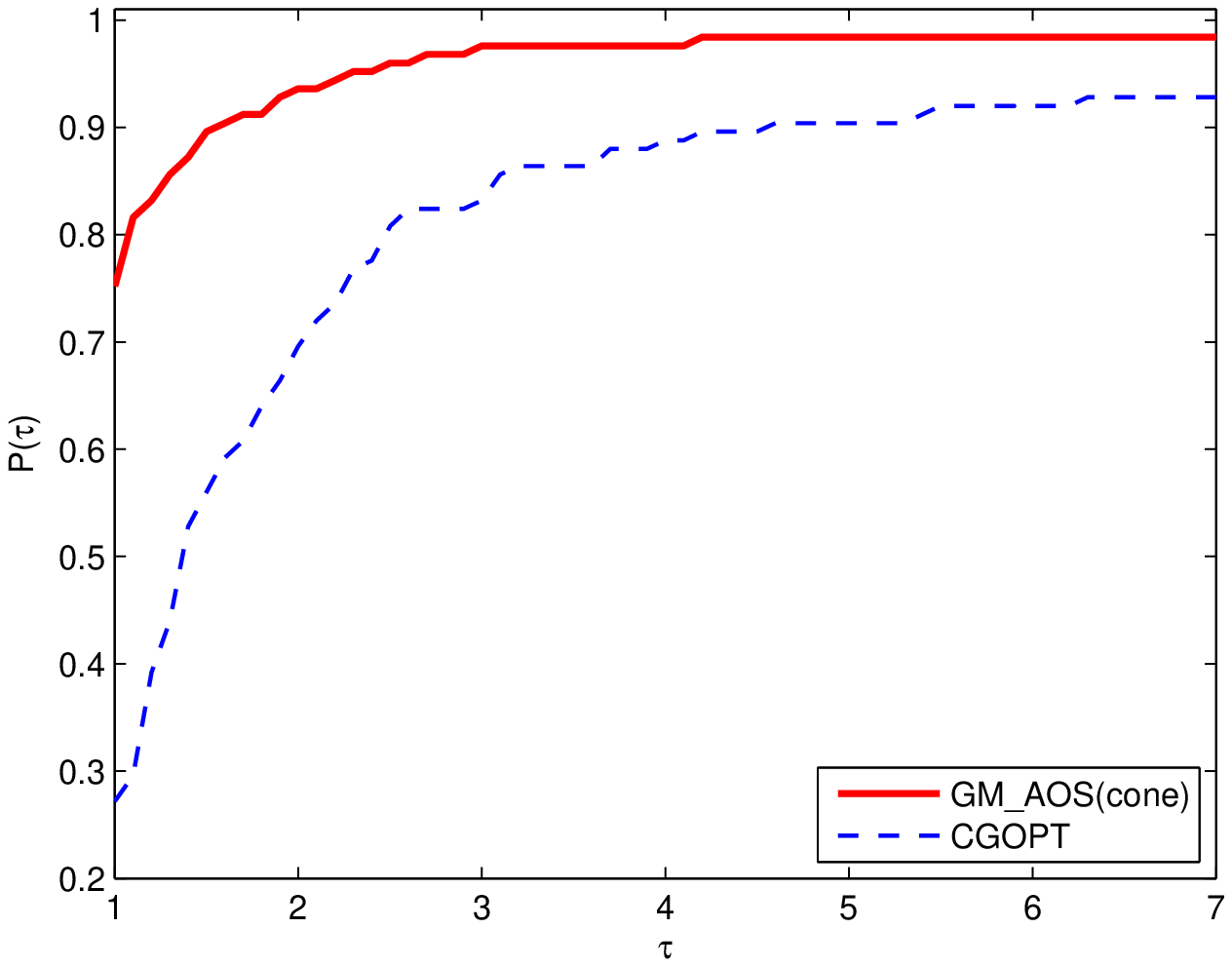}
                 		\caption{   $ N_{f} $ (144pCUTEr)}\label{fig:5-13}
                 	\end{minipage}%
                 	\begin{minipage}[t]{0.5\linewidth}
                 		\includegraphics[width=8cm,height=5.5cm]{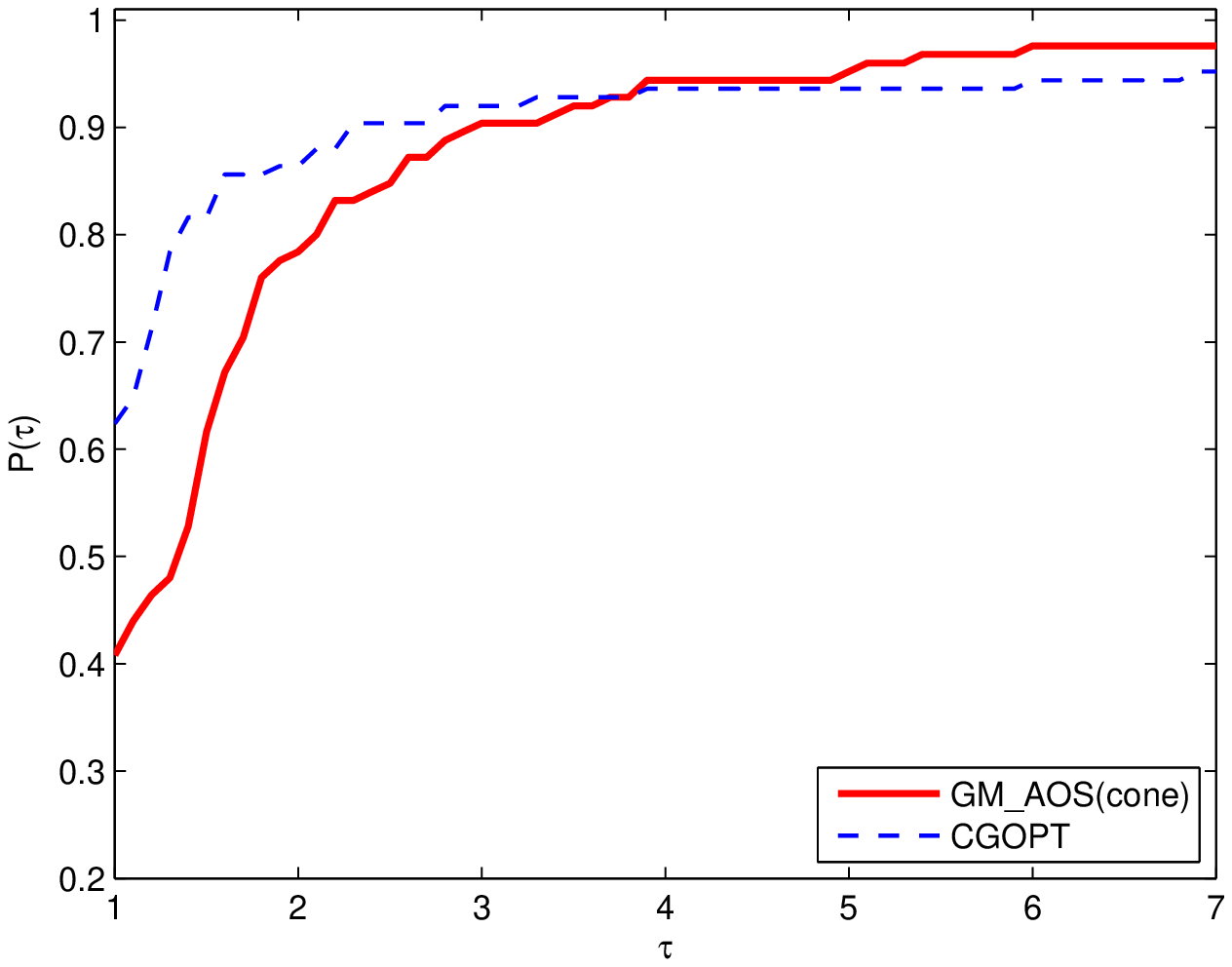}
                 		\caption{   $ N_{g} $ (144pCUTEr)}
                 		\label{fig:5-14}
                 	\end{minipage}%
                 \end{figure*}  	
                 
                 \begin{figure*}[htbp]
                 	\begin{minipage}[t]{0.5\linewidth}
                 		\includegraphics[width=8cm,height=5.5cm]{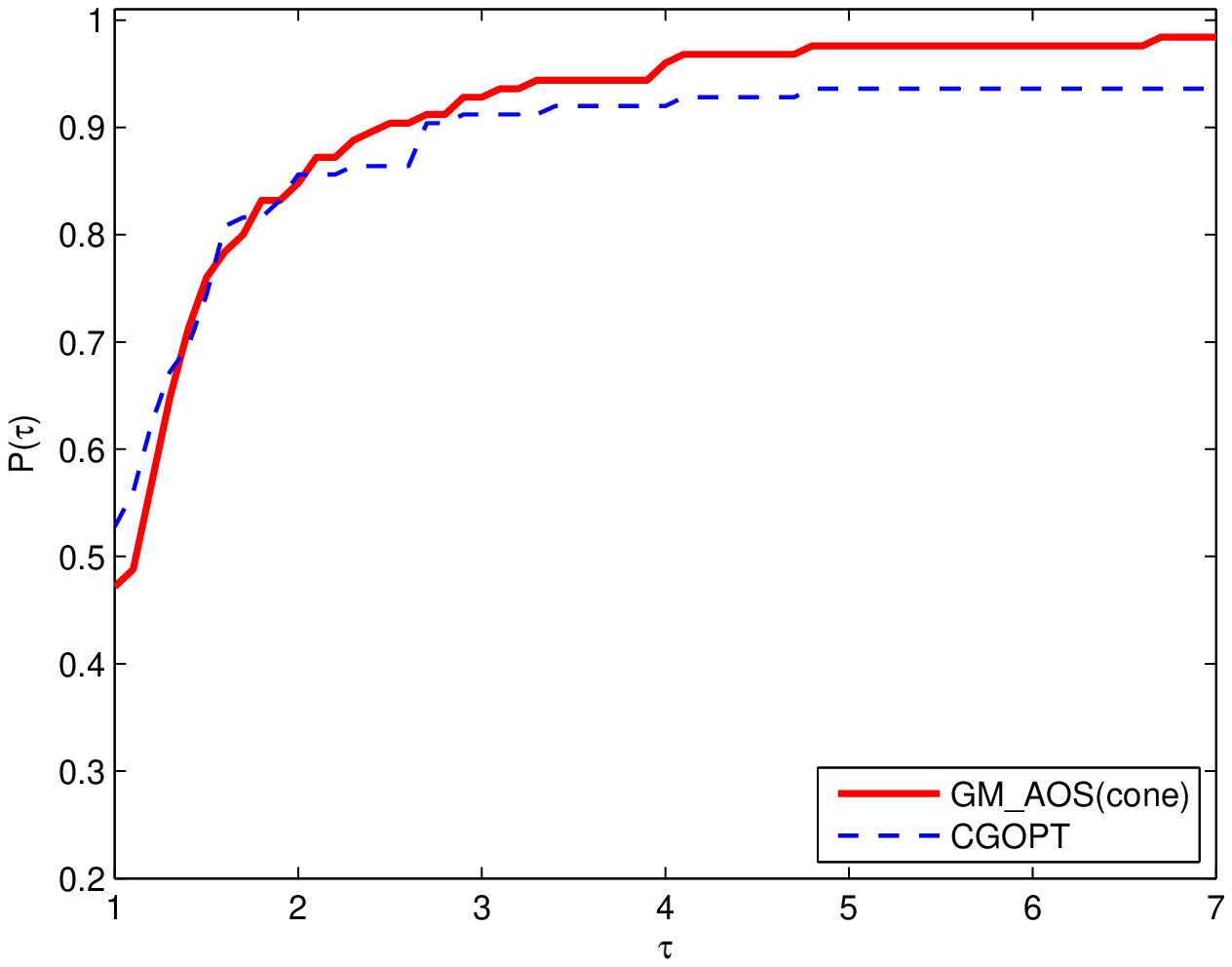}
                 		\caption{$ N_f+3N_g $(144pCUTEr)}\label{fig:5-15}
                 	\end{minipage}%
                 	\begin{minipage}[t]{0.5\linewidth}
                 		\includegraphics[width=8cm,height=5.5cm]{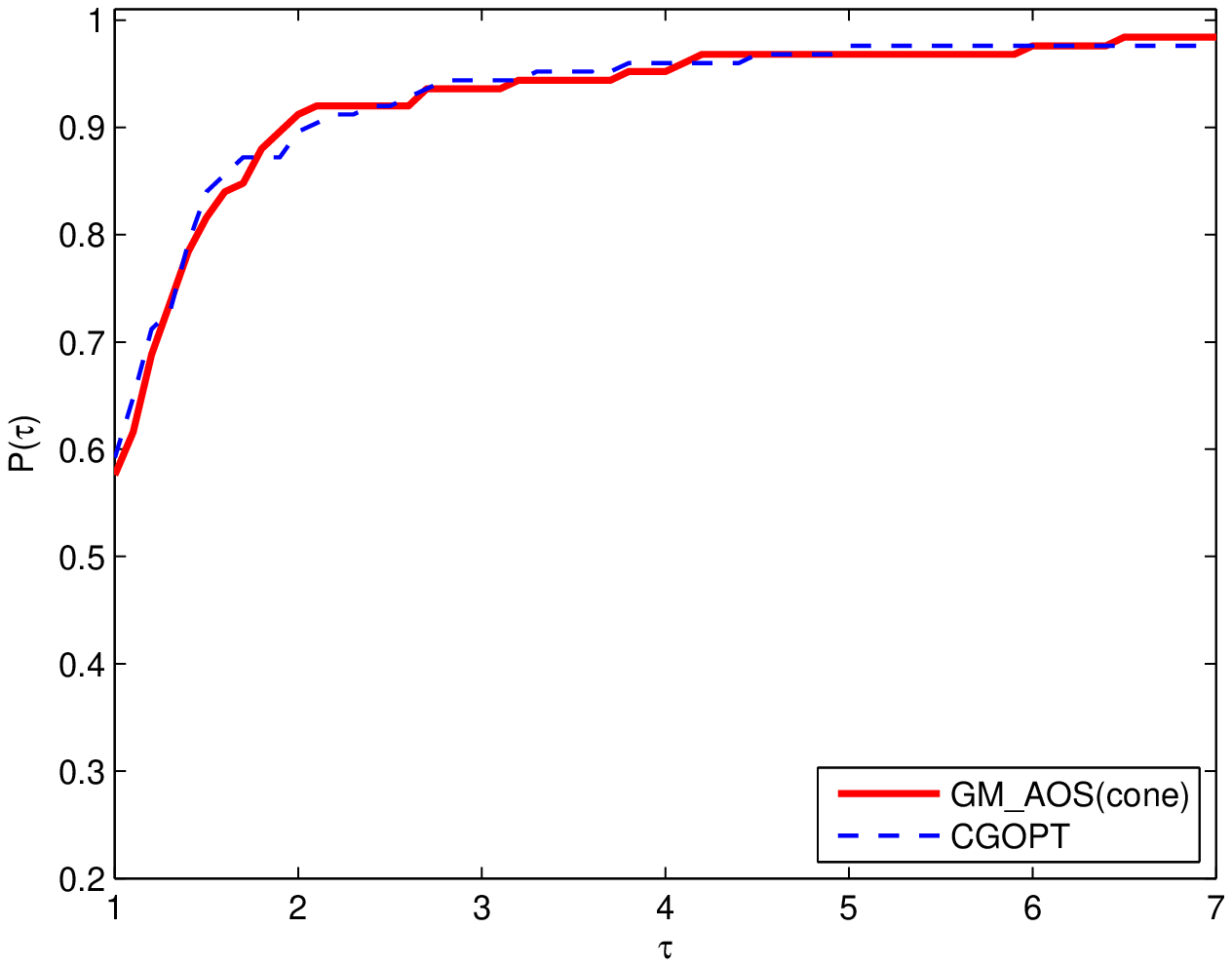}
                 		\caption{   $ T_{CPU} $ (144pCUTEr)}
                 		\label{fig:5-16}
                 	\end{minipage}%
                 \end{figure*}

In  the third group of numerical experiments,   we use the collect sets 80pAndr and 144pCUTEr to  compare  the performance of    GM$ \_$AOS    with that of  CGOPT.  For 80pAndr, GM$ \_ $AOS (cone) successfully solves all 80 problems, while  CGOPT  successfully solves 79 problems.
 Figs. \ref{fig:5-9}-\ref{fig:5-12} plot the performance profiles of  GM$ \_ $AOS (cone) and CGOPT  for 80pAndr in term of $ N_t $,  $ N_f $, $ N_g $ and $ T_{CPU} $.  As shown in Figs. \ref{fig:5-9}-\ref{fig:5-12}, we observe that GM$ \_ $AOS (cone)  is considerably superior to CGOPT for 80pAndr.  For 144pCUTEr, GM$ \_ $AOS (cone) successfully solves 134 problems, while  CGOPT  successfully solves 133 problems.  Figs. \ref{fig:5-13}-\ref{fig:5-16} plot the performance profiles of  GM$ \_ $AOS (cone) and CGOPT  for 144pCUTEr in term of $ N_f $,  $ N_g $, $ N_f + 3 N_g $ \cite{ZhangHager2006Algorithm}  and $ T_{CPU} $.  As shown in Fig. \ref{fig:5-13},  GM$ \_ $AOS (cone)  is much superior to CGOPT in term of $ N_f $, since GM$ \_ $AOS (cone) solves about 78$ \% $ problems with the least function evaluation, while the percentage of CGOPT   is about 29$ \% $ for 80pAndr. Fig. \ref{fig:5-14} indicates that GM$ \_ $AOS (cone) is inferior to CGOPT in term of $ N_g $, while Fig. \ref{fig:5-15} shows that GM$ \_ $AOS (cone) performs a little better than CGOPT in term of $ N_f + 3 N_g $.  We observe from Fig. \ref{fig:5-16} that GM$ \_ $AOS (cone) is as fast as CGOPT.   It indicates GM$ \_ $AOS (cone)  is much superior to CGOPT for 80pAdr and is comparable  to  CGOPT for  144pCUTEr.

                 \begin{figure*}[htbp]
                 	\begin{minipage}[t]{0.5\linewidth}
                 		\includegraphics[width=8cm,height=5.5cm]{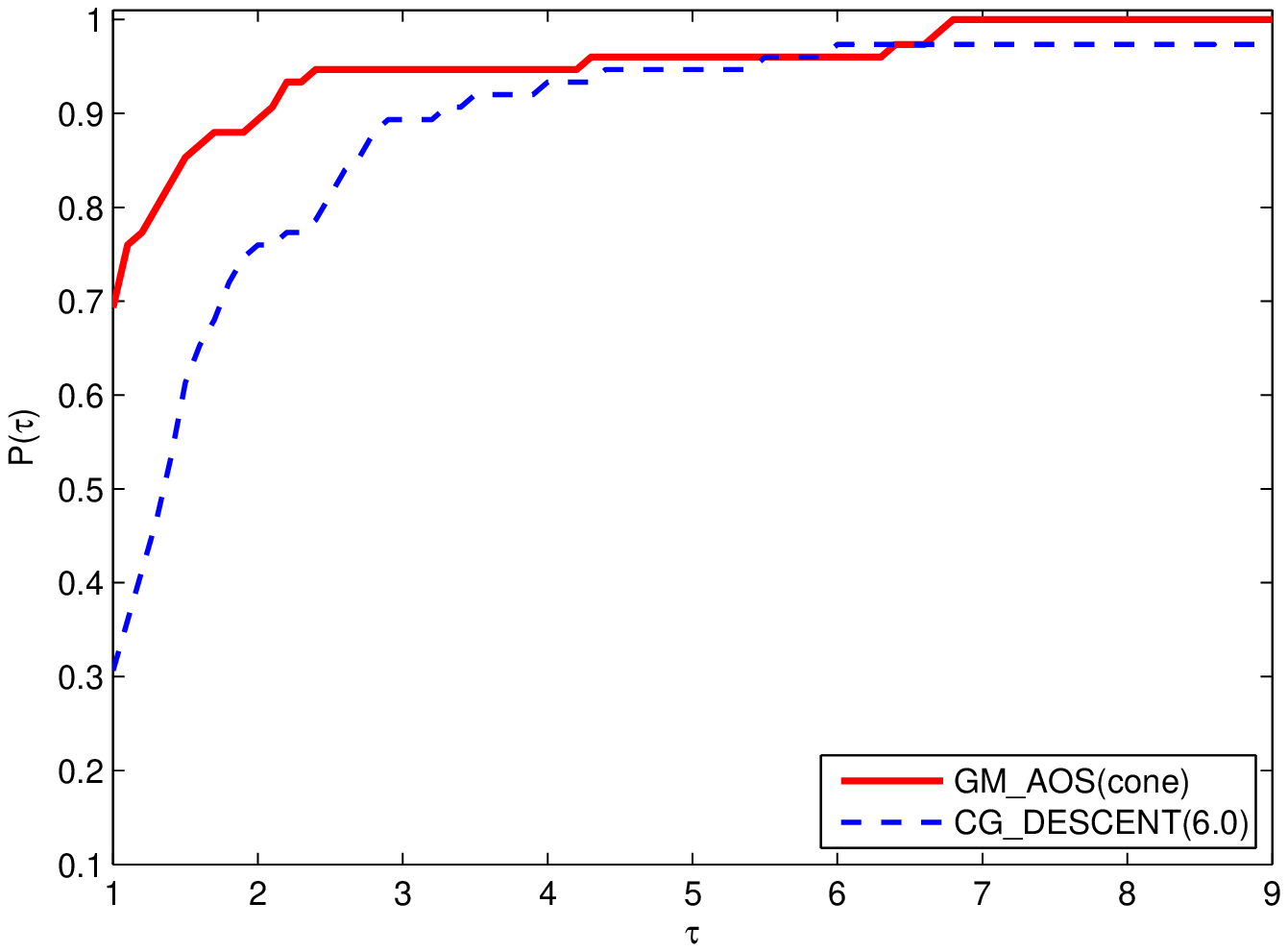}
                 		\caption{   $ N_{f} $            (80pAndr)}\label{fig:5-17}
                 	\end{minipage}%
                 	\begin{minipage}[t]{0.5\linewidth}
                 		\includegraphics[width=8cm,height=5.5cm]{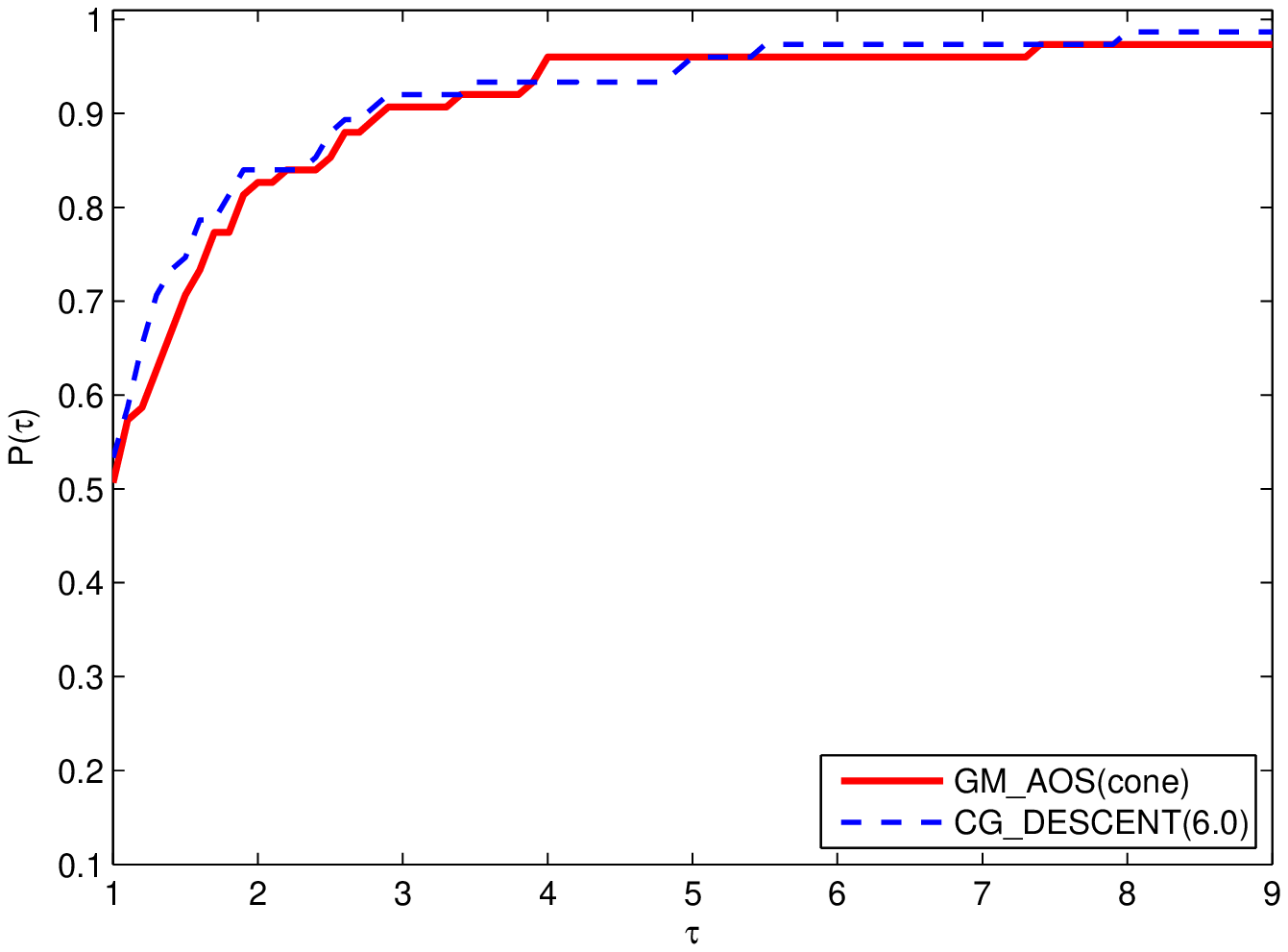}
                 		\caption{   $ N_{g} $ (80pAndr)}
                 		\label{fig:5-18}
                 	\end{minipage}%
                 \end{figure*}  	
                 \begin{figure*}[htbp]
                 	\begin{minipage}[t]{0.5\linewidth}
                 		\includegraphics[width=8cm,height=5.5cm]{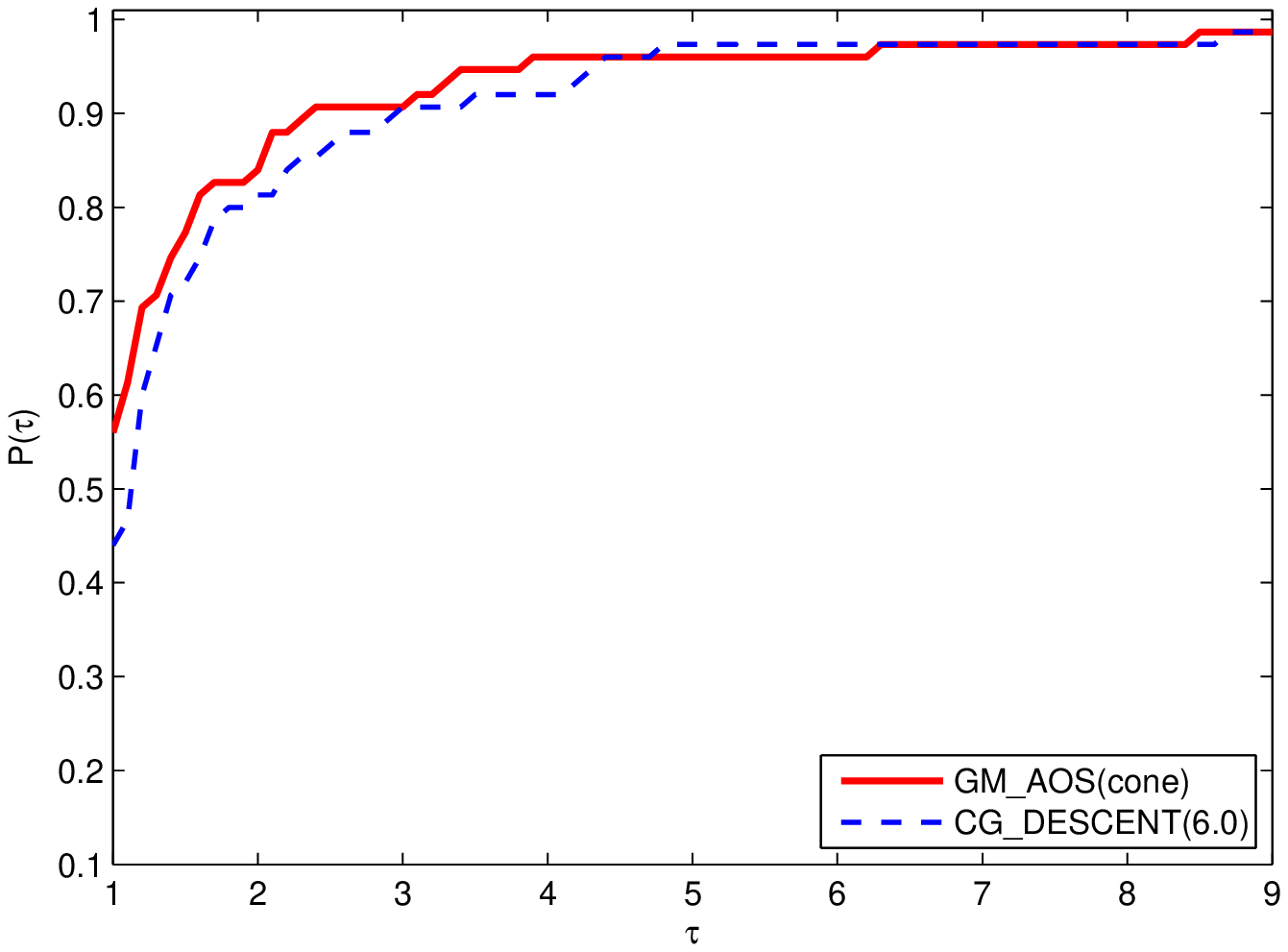}
                 		\caption{ $ N_f+3N_g $(80pAndr)}\label{fig:5-19}
                 	\end{minipage}%
                 	\begin{minipage}[t]{0.5\linewidth}
                 		\includegraphics[width=8cm,height=5.5cm]{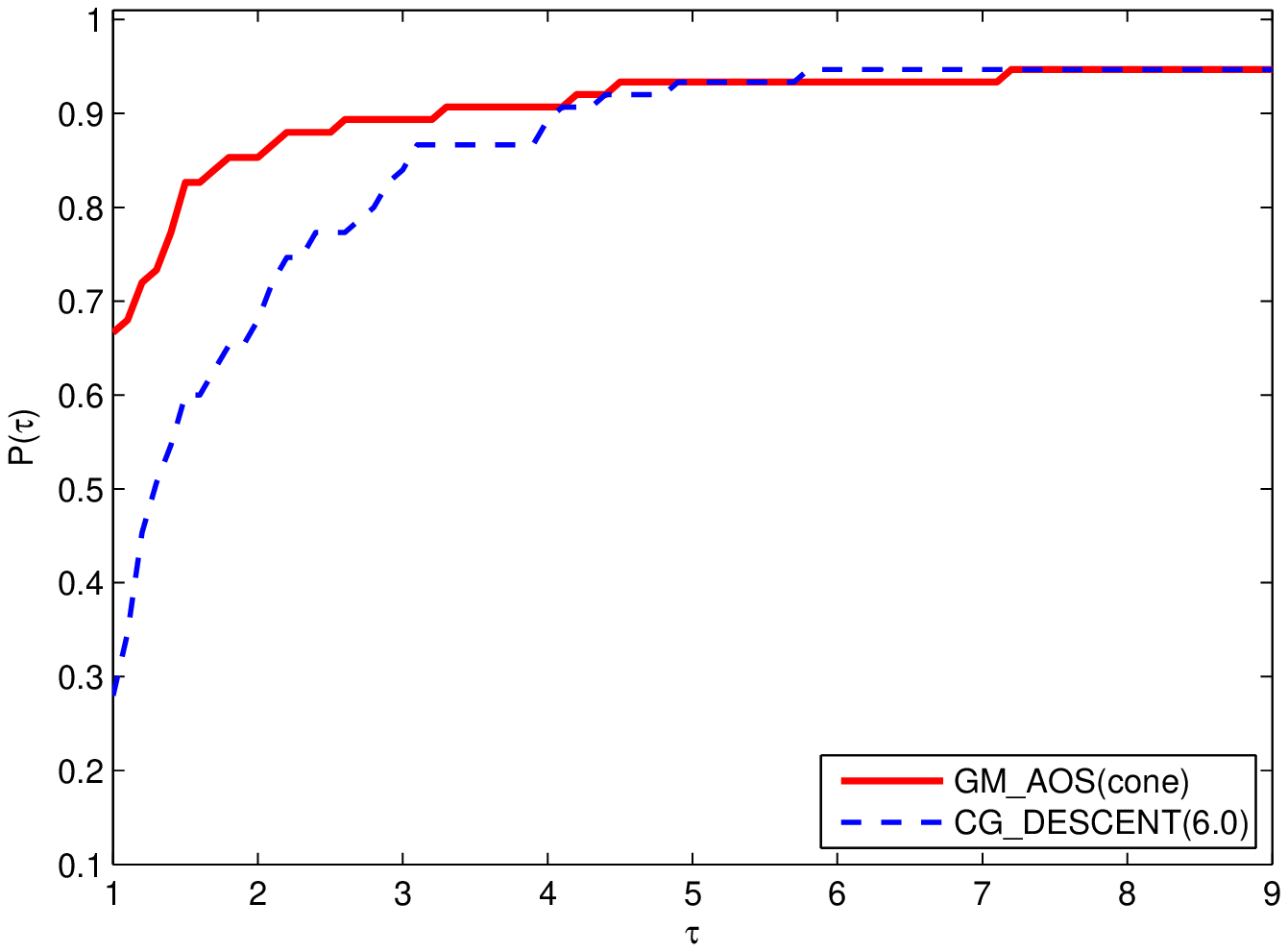}
                 		\caption{$ T_{CPU} $(80pAndr)}
                 		\label{fig:5-20}
                 	\end{minipage}%
                 \end{figure*}

In  the fourth group of numerical experiments,   we  use the collect set 80pAndr to  compare  the performance of    GM$ \_$AOS    with that of  CG$ \_ $DESCENT (6.0), which is the limited memory conjugate gradient software package,  and then use the collect set 144pCUTEr to  compare  the performance of    GM$ \_$AOS    with that of  CG$ \_ $DESCENT (5.0).  For 80pAndr, GM$ \_ $AOS (cone)  successfully solves all 80 problems, while  CG$ \_ $DESCENT (6.0)  successfully solves 75 problems.
Figs. \ref{fig:5-17}-\ref{fig:5-20} plot the performance profiles of  GM$ \_ $AOS (cone) and CG$ \_ $DESCENT   (6.0)  for 80pAndr in term of $ N_f $,  $ N_g $, $  N_f+3 N_g $ \cite{ZhangHager2006Algorithm} and $ T_{CPU} $.  As shown in Figs.\ref{fig:5-17}-\ref{fig:5-18}, we observe that GM$ \_ $AOS (cone)  is considerably superior to CG$ \_ $DESCENT (6.0)  in term of $ N_f $ but is  a little inferior to CG$ \_ $DESCENT (6.0) in term of $ N_g $. Fig. \ref{fig:5-19} shows that GM$ \_ $AOS (cone) performs better than CG$ \_ $DESCENT (6.0) in term of $N_f+3 N_g$. We observe from Fig. \ref{fig:5-20} that GM$ \_ $AOS (cone) is   faster than CG$ \_ $DESCENT (6.0).  For 144pCUTEr, GM$ \_ $AOS (cone) successfully solves 134 problems, while  CG$ \_ $DESCENT (5.0)  successfully solves 142 problems.  Figs. \ref{fig:5-21}-\ref{fig:5-24} plot the performance profiles of  GM$ \_ $AOS (cone) and CG$ \_ $DESCENT (5.0)  for 144pCUTEr in term of $ N_f $,  $ N_g $, $ N_f + 3 N_g $ and $ T_{CPU} $.  As shown in Fig. \ref{fig:5-21},  GM$ \_ $AOS (cone)  is much superior to CGOPT in term of $ N_f $, since GM$ \_ $AOS (cone) solves about 65$ \% $ problems with the least function evaluations, while  the percentage of CG$ \_ $DESCENT (5.0)  is about 39$ \% $.   Fig. \ref{fig:5-22} indicates that GM$ \_ $AOS (cone) is inferior to CG$ \_ $DESCENT (5.0) in term of $ N_g $, while Fig. \ref{fig:5-23} shows that GM$ \_ $AOS (cone) is comparable to CG$ \_ $DESCENT (5.0) in term of $ N_f+3N_g $.  We observe from Fig. \ref{fig:5-24} that GM$ \_ $AOS (cone) is as fast as CG$ \_ $DESCENT (5.0).   It indicates GM$ \_ $AOS (cone)  is   superior to CG$ \_ $DESCENT (6.0) for 80pAdr and is comparable  to  CG$ \_ $DESCENT (5.0) for  144pCUTEr.
                 \begin{figure*}[htbp]
                 	\begin{minipage}[t]{0.5\linewidth}
                 		\includegraphics[width=8cm,height=5.5cm]{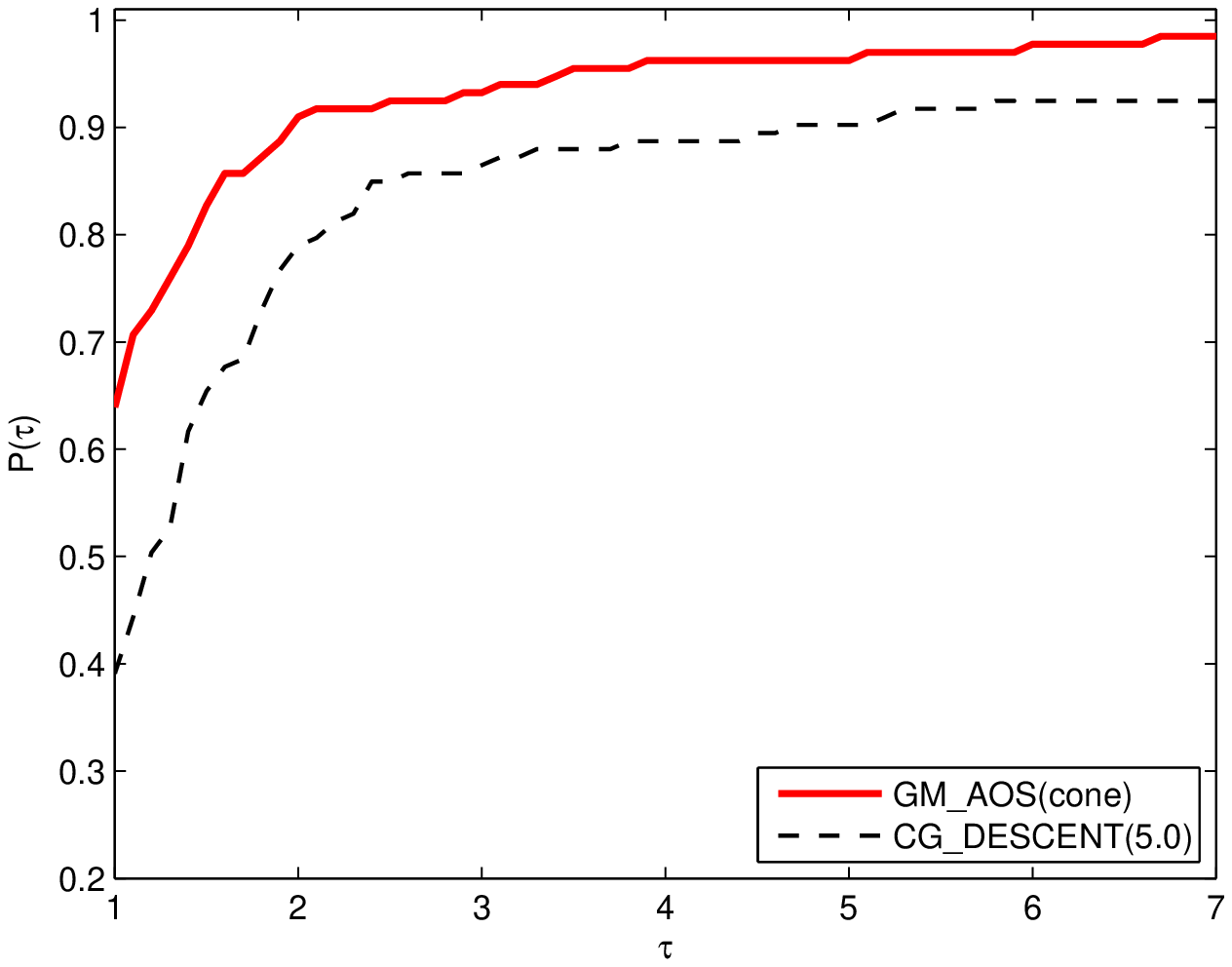}
                 		\caption{  $ N_f $ (144pCUTEr)}\label{fig:5-21}
                 	\end{minipage}%
                 	\begin{minipage}[t]{0.5\linewidth}
                 		\includegraphics[width=8cm,height=5.5cm]{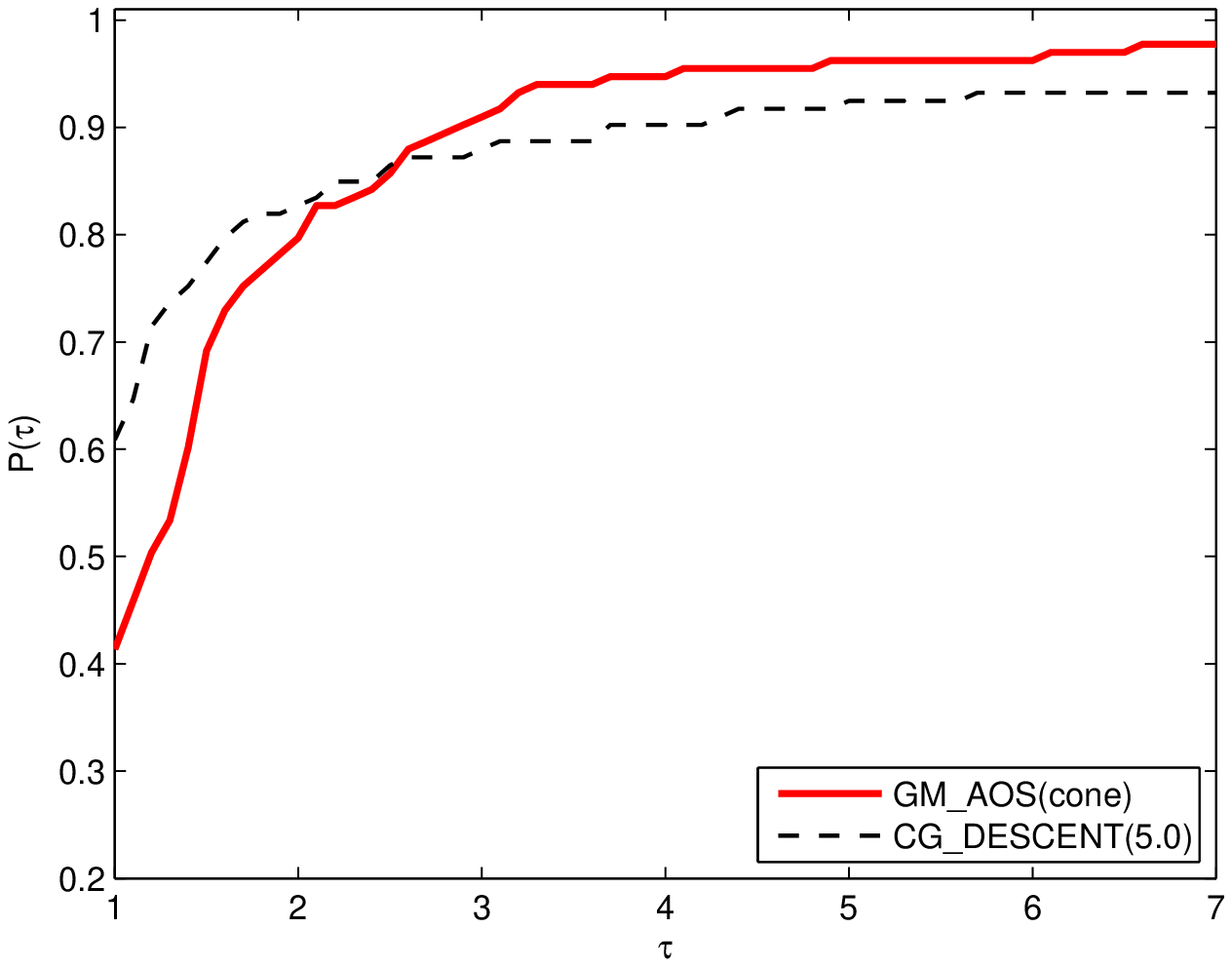}
                 		\caption{  $ N_g $ (144pCUTEr)}
                 		\label{fig:5-22}
                 	\end{minipage}%
                 \end{figure*}
                 
                 \begin{figure*}[htbp]
                 	\begin{minipage}[t]{0.5\linewidth}
                 		\includegraphics[width=8cm,height=5.5cm]{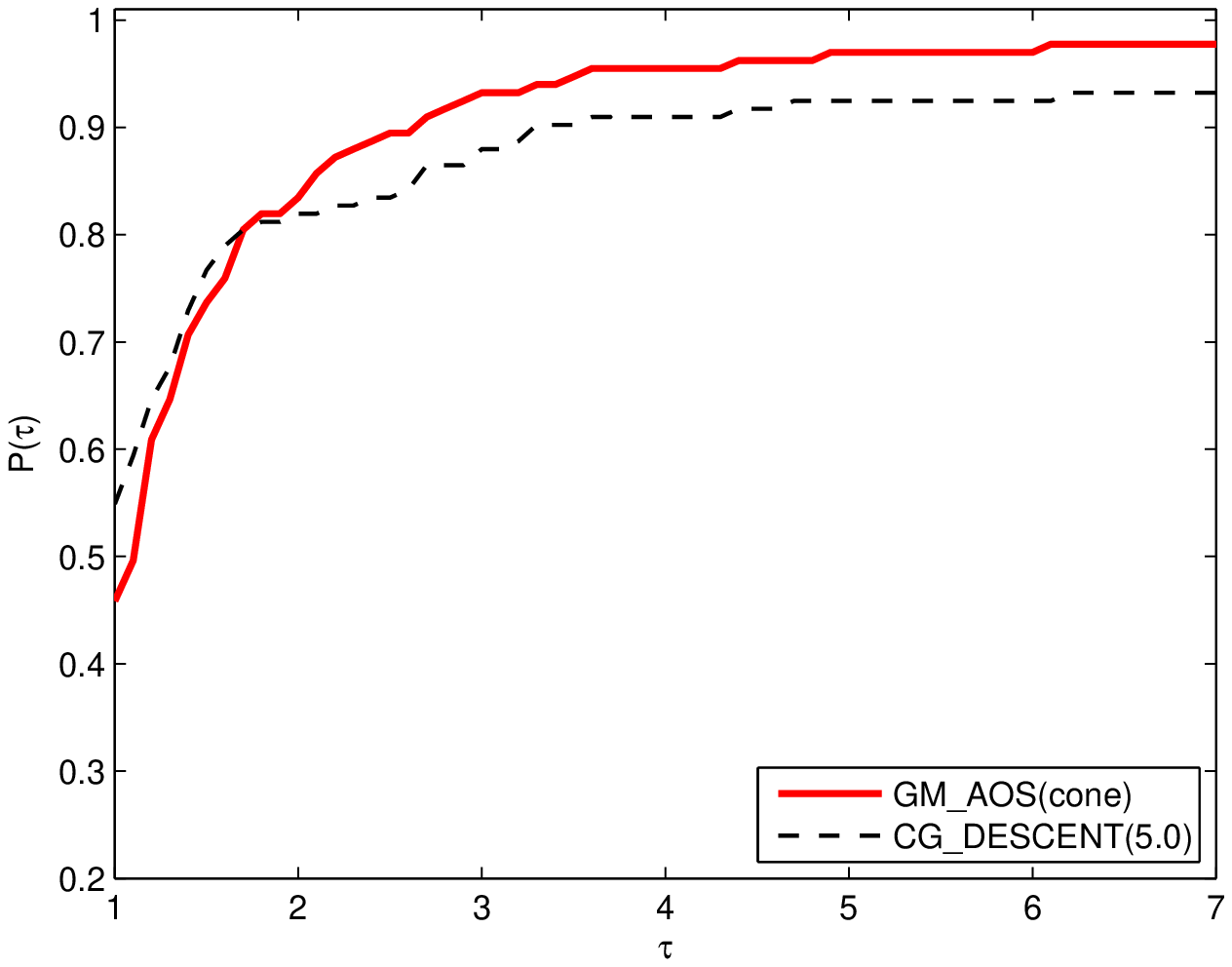}
                 		\caption{ $ N_f+3N_g  $(144pCUTEr)}\label{fig:5-23}
                 	\end{minipage}%
                 	\begin{minipage}[t]{0.5\linewidth}
                 		\includegraphics[width=8cm,height=5.5cm]{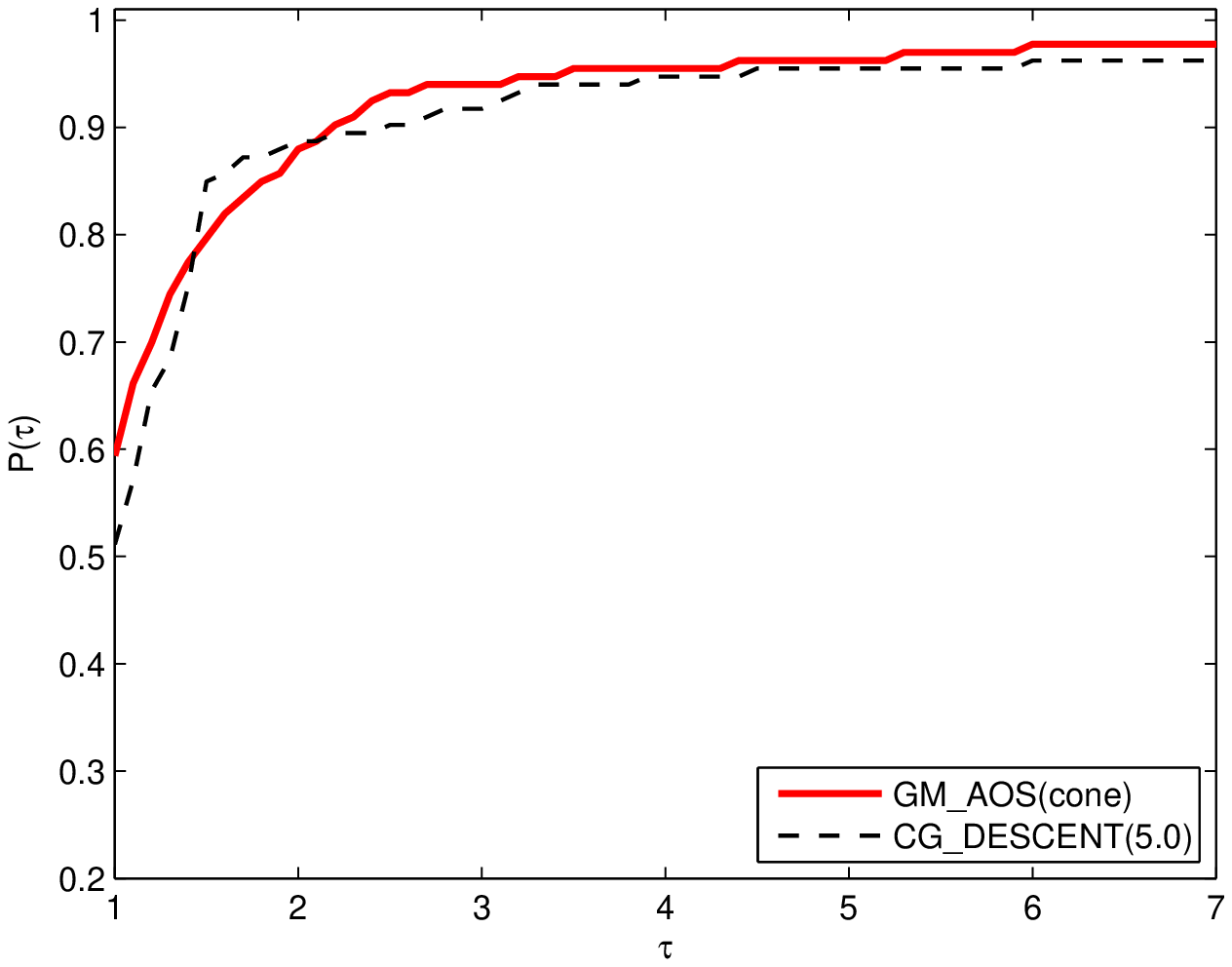}
                 		\caption{   $ T_{CPU} $ (144pCUTEr)}
                 		\label{fig:5-24}
                 	\end{minipage}%
                 \end{figure*}

\section{Conclusions and Discussions}
\indent {In this paper, we present an improved   gradient method with approximately optimal stepsize based on conic model (GM$ \_ $AOS (cone)). In GM$ \_ $AOS (cone),  some approximation models including the conic model and some quadratic models are exploited to generate   approximately optimal stepsizes for gradient method. It is noted that the main difference between the proposed method and the gradient method with approximately opitmal stepsize based on conic model \cite{Liu2018GMAOScone} lies that  the proposed method uses the stepsize (3.1) as the initial stepsize at the first iteration, while the gradient method \cite{Liu2018GMAOScone} takes $1/\left\| g_0\right\|_{\infty}   $ as  the initial stepsize. In addition, more numerical experiments with two group collect sets 80pAdr and 144pCUTEr  are conducted to examine the effectiveness of the proposed method.  Numerical results indicate that GM$ \_ $AOS (cone)  is    superior to the SBB4 method and the BB method,   performs better than CGOPT \cite{Dai2013CGOPT} for   80pAdr and is comparable to CGOPT for  144pCUTEr, and  is   superior to the limited memory conjugate gradient software package CG$ \_ $DESCENT (6.0) \cite{HagerZhang2013The} for   80pAdr 
 \begin{table}[!htbp]
 	\centering
 	\renewcommand{\arraystretch}{0.7}
 	\caption{ The test functions in collect set 80pAndr  mainly from   \cite{Andrei2008An}}
 	\begin{tabular}{{l l }}
 		\toprule		\label{tb:5-1}
 		Name &	 	Name  	\\
 		\midrule
 		Freudenstein and Roth   FREUROTH (CUTE)     	&	EG2 (CUTE)	\\
 		Extended Trigonometric ET1	&	EDENSCH (CUTE)	\\
 		Extended Rosenbrock   SROSENBR (CUTE)	&	Broyden Pentadiagonal (CUTE)	\\
 		Extended White and Holst	&	Almost Perturbed Quadratic 	\\
 		Extended Beale  BEALE (CUTE)	&	Almost Perturbed Quartic	\\
 		Extended Penalty	&	FLETCHCR (CUTE)	\\
 		Perturbed Quadratic	&	ENGVAL1 (CUTE)	\\
 		Raydan 1	&	DENSCHNA (CUTE)	\\
 		Raydan 2	&	DENSCHNB (CUTE)	\\
 		TR-SUMM     	&	DENSCHNC (CUTE)	\\
 		Diagonal 1  	&	DENSCHNF (CUTE)	\\
 		Diagonal 2  	&	SINQUAD  (CUTE)	\\
 		Hager   	&	HIMMELBG (CUTE)	\\
 		Generalized Tridiagonal 1	&	HIMMELBH (CUTE)	\\
 		Extended Tridiagonal 1	&	DIXON3DQ (CUTE)	\\
 		Extended Three Expo Terms	&	BIGGSB1 (CUTE)	\\
 		Generalized Tridiagonal 2	&	Perturbed Quadratic	\\
 		Diagonal 3 (1c1c)	&	GENROSNB (CUTE)	\\
 		Diagonal Full Borded	&	QP1 Extended Quadratic Penalty	\\
 		Extended Himmelblau  HIMMELBC (CUTE) 	&	QP2 Extended Quadratic Penalty 	\\
 		Extended Powell	&	Tridiagonal TS1 	\\
 		Tridiagonal Double Borded Arrow Up            	&	Tridiagonal TS2	\\
 		Extended PSC1	&	Tridiagonal TS3	\\
 		Extended Block-Diagonal BD1	&	Extended Trigonometric ET2	\\
 		Extended Maratos	&	QP3 Extended Quadratic Penalty 	\\
 		Full Hessian FH1	&	EG1	\\
 		Extended Cliff	&	GENROSEN-2	\\
 		Quadratic Diagonal Perturbed	&	PRODsin	\\
 		Full Hessian FH2	&	PROD1 (m=n)	\\
 		Full Hessian FH3	&	PRODcos                 	\\
 		Tridiagonal Double Borded - NONDQUAR	&	PROD2 (m=1)	\\
 		Tridiagonal White and Holst (c=4)	&	ARGLINB (m=5)	\\
 		Diagonal Double Borded Arrow Up	&	DIXMAANA (CUTE)	\\
 		TRIDIA (CUTE)	&	DIXMAANB (CUTE)	\\
 		ARWHEAD  (CUTE) 	&	DIXMAANC (CUTE)	\\
 		NONDIA  (CUTE)	&	DIXMAAND (CUTE)	\\
 		Extended WOODS (CUTE)	&	DIXMAANL (CUTE)	\\
 		Extended Hiebert	&	VARDIM (CUTE)	\\
 		BDQRTIC (CUTE)	&	DIAG-AUP1 	\\
 		DQDRTIC (CUTE)		&	ENGVAL8                                 	\\
 		\bottomrule	
 	\end{tabular}
 \end{table}

 \begin{table}[!htbp]
 	\centering
 	\renewcommand{\arraystretch}{0.7}
 	\caption{ The test functions in collect set 144pCUTEr  from CUTEr library \cite{Gould2001CUTEr}}
 	\begin{tabular}{{lccc cccc c}}
 		\toprule	\label{tb:5-2}
 		Name &	Dimension	&	Name	&	Dimension	&	Name	&	Dimension 	\\
 		\midrule
 		AKIVA	&	2	&	EDENSCH	&	2000	&	NONDQUAR	&	5000	\\
 		ALLINITU	&	4	&	EG2	&	1000	&	OSBORNEA	&	5	\\
 		ARGLINA	&	200	&	EIGENALS	&	2550	&	OSBORNEB	&	11	\\
 		ARGLINB	&	200	&	EIGENBLS	&	2550	&	OSCIPATH	&	10	\\
 		ARWHEAD	&	5000	&	EIGENCLS	&	2652	&	PALMER1C	&	8	\\
 		BARD	&	3	&	ENGVAL1	&	5000	&	PALMER1D	&	7	\\
 		BDQRTIC	&	5000	&	ENGVAL2	&	3	&	PALMER2C	&	8	\\
 		BEALE	&	2	&	ERRINROS	&	50	&	PALMER3C	&	8	\\
 		BIGGS6	&	6	&	EXPFIT	&	2	&	PALMER4C	&	8	\\
 		BOX3	&	3	&	EXTROSNB	&	1000	&	PALMER5C	&	6	\\
 		BOX	&	10000	&	FLETCBV2	&	5000	&	PALMER6C	&	8	\\
 		BRKMCC	&	2	&	FLETCHCR	&	1000	&	PALMER7C	&	8	\\
 		BROWNAL	&	200	&	FMINSRF2	&	5625	&	PALMER8C	&	8	\\
 		BROWNBS	&	2	&	FMINSURF	&	5625	&	PARKCH	&	15	\\
 		BROWNDEN	&	4	&	FREUROTH	&	5000	&	PENALTY1	&	1000	\\
 		BROYDN7D	&	5000	&	GENHUMPS	&	5000	&	PENALTY2	&	200	\\
 		BRYBND	&	5000	&	GENROSE	&	500	&	PENALTY3	&	200	\\
 		CHAINWOO	&	4000	&	GROWTHLS	&	3	&	POWELLSG	&	5000	\\
 		CHNROSNB	&	50	&	GULF	&	3	&	POWER	&	10000	\\
 		CLIFF	&	2	&	HAIRY	&	2	&	QUARTC	&	5000	\\
 		COSINE	&	10000	&	HATFLDD	&	3	&	ROSENBR	&	2	\\
 		CRAGGLVY	&	5000	&	HATFLDE	&	3	&	S308	&	2	\\
 		CUBE	&	2	&	HATFLDFL	&	3	&	SCHMVETT	&	5000	\\
 		CURLY10	&	10000	&	HEART6LS	&	6	&	SENSORS	&	100	\\
 		CURLY20	&	10000	&	HEART8LS	&	8	&	SINEVAL	&	2	\\
 		CURLY30	&	10000	&	HELIX	&	3	&	SINQUAD	&	5000	\\
 		DECONVU	&	63	&	HIELOW	&	3	&	SISSER	&	2	\\
 		DENSCHNA	&	2	&	HILBERTA	&	2	&	SNAIL	&	2	\\
 		DENSCHNB	&	2	&	HILBERTB	&	10	&	SPARSINE	&	5000	\\
 		DENSCHNC	&	2	&	HIMMELBB	&	2	&	SPARSQUR	&	10000	\\
 		DENSCHND	&	3	&	HIMMELBF	&	4	&	SPMSRTLS	&	4999	\\
 		DENSCHNE	&	3	&	HIMMELBG	&	2	&	SROSENBR	&	5000	\\
 		DENSCHNF	&	2	&	HIMMELBH	&	2	&	STRATEC	&	10	\\
 		DIXMAANA	&	3000	&	HUMPS	&	2	&	TESTQUAD	&	5000	\\
 		DIXMAANB	&	3000	&	JENSMP	&	2	&	TOINTGOR	&	50	\\
 		DIXMAANC	&	3000	&	JIMACK	&	3549	&	TOINTGSS	&	5000	\\
 		DIXMAAND	&	3000	&	KOWOSB	&	4	&	TOINTPSP	&	50	\\
 		DIXMAANE	&	3000	&	LIARWHD	&	5000	&	TOINTQOR	&	50	\\
 		DIXMAANF	&	3000	&	LOGHAIRY	&	2	&	TQUARTIC	&	5000	\\
 		DIXMAANG	&	3000	&	MANCINO	&	100	&	TRIDIA	&	5000	\\
 		DIXMAANH	&	3000	&	MARATOSB	&	2	&	VARDIM	&	200	\\
 		DIXMAANI	&	3000	&	MEXHAT	&	2	&	VAREIGVL	&	50	\\
 		DIXMAANJ	&	3000	&	MOREBV	&	5000	&	VIBRBEAM	&	8	\\
 		DIXMAANK	&	15	&	MSQRTALS	&	1024	&	WATSON	&	12	\\
 		DIXMAANL	&	3000	&	MSQRTBLS	&	1024	&	WOODS	&	4000	\\
 		DIXON3DQ	&	10000	&	NCB20B	&	5000	&	YFITU	&	3	\\
 		DJTL	&	2	&	NCB20	&	5010	&	ZANGWIL2	&	2	\\
 		DQDRTIC	&	5000	&	NONCVXU2	&	5000	&		&		\\
 		DQRTIC	&	5000	&	NONDIA	&	5000	&		&		\\	 	
 		\bottomrule	
 	\end{tabular}
 \end{table}

\noindent	and is comparable to CG$ \_ $DESCENT (5.0) \cite{Hager2005A} for   144pCUTEr. As far as we know, GM$ \_ $AOS (cone) is the most efficient gradient method for general unconstrained optimization so far.
	
	Given that the search direction $-g_k$     has low storage and can be easily computed,   the  nonmonotone Armijo  line search used   can be easily implemented and the numerical effect is   surprising,  the gradient methods with approximately optimal stepsizes will be   strong candidates  for  large scale  unconstrained optimization. And the following problems are very interesting:(1)		What is the best gradient method with approximately optimal stepsize (GM$ \_ $AOS) ?	(2) Can the gradient method with approximately optimal stepsize perform better than CG$ \_ $DESCENT (5.3) for CUTEr library?

\begin{acknowledgements}  We  would like to thank Professor  Dai, Y. H. and Dr. Kou caixia for their  C  code of CGOPT, and thank Professors Hager and Zhang, H. C. for their C code of CG$ \_ $DESCENT. This research is supported by    Guangxi Natural Science Foundation (No.2018GXNSFBA281180). \newline
\end{acknowledgements}

\end{document}